\documentclass[a4paper, 12pt]{article}
\usepackage{times}
\usepackage[latin1]{inputenc}
\usepackage[T1]{fontenc}
\usepackage{amssymb}
\usepackage{amsmath}
\usepackage{amsthm}
\usepackage{latexsym}
\usepackage[dvips]{graphicx}
\usepackage[margin=3cm]{geometry}
\def\C{{\mathbb C}}
\def\Z{{\mathbb Z}}
\def\N{{\mathbb N}}
\def\R{{\mathbb R}}
\def\T{{\mathbb T}} 
\def\e{{\epsilon}}
\def\g{{\gamma}}
\def\s{{\sigma}}
\def\r{{\rho}}
\def\o{{\omega}}
\def\a{{\alpha}}
\def\d{{\delta}}
\def\l{{\lambda}}
\def\p{{\prime}}
\def\for{{\,\,\, \forall \,\,}}

\def\B{ {\mathcal{B} } }
\def\A{ {\mathcal{A} } }
\def\K{ {\mathcal{K} } }

\def\mS{ {\mathcal{S} } }
\def\Y{ {\mathcal{Y} } }

\def\Im{\,{\rm Im}\,}
\def\Re{\,{\rm Re}\,}

\def\meas{\,{\rm meas}}

\def\Russ{R\"{u}{\ss}mann}
\def\+R{+_{_{ \!\! \R}}}

\def\Kep{{\mbox{\scriptsize{Kep}}}}

\def\M{{\mathcal{M}}}
\def\vs{{\varsigma}}

\def\diag{\,\mbox{diag}\,}

\newtheorem{teo}{Theorem}
\newtheorem{pro}{Proposition}
\newtheorem{lem}{Lemma}
\newtheorem{defin}{Definition}
\newtheorem{oss}{Remark}

\def\nl{\vglue0.2truecm\noindent}

\begin{document}

\title{Analytic Lagrangian tori for the planetary\\ many--body problem\thanks{{\bf Acknowledgments.} We are indebted with 
Jacques F\'{e}joz for many enlightening discussions. This work was partially supported
by  the Italian MIUR project ``Metodi variazionali e equazioni differenziali nonlineari''. }}

\author{Luigi Chierchia and Fabio Pusateri}

\author{ 
Luigi Chierchia
\and
Fabio Pusateri
}

\date{}

\maketitle

\abstract\noindent{\footnotesize
In 2004 J. F\'{e}joz \cite{fejoz}, completing investigations of M. Herman's \cite{herman},
gave a complete proof  of ``Arnold's Theorem'' \cite{arnold} on the  planetary many--body problem, establishing, in particular, 
the existence of a positive measure set of  smooth ($C^\infty$) Lagrangian invariant tori  for the planetary many--body problem.
Here, using \Russ's 2001 KAM theory \cite{russmann}, we prove the above result in the real--analytic class.}

\tableofcontents

\section{Introduction}
\label{risultati}

\phantom{A}

The planetary many--body problem consists in studying the evolution of $(1+n)$ bodies (point masses), subject only to the mutual gravitational attraction, in the case one of bodies (the ``Sun'') has mass $m_0$ considerably larger than the masses $m_i$  of the remaining $n$ bodies (the ``planets''). The Newtonian evolution equations  
for such problem (in suitable units) are given by
\begin{equation}\label{Newton}
m_j \
\ddot q_j = \sum_{k\neq j} m_j m_k\ \frac{q_k-q_j}{|q_k-q_j|^3}\ ,\qquad j=0,1,...,n\ ,
\end{equation}
where $q_j=q_j(t)\in\R^3$ denotes the position at time $t$ of the $j^{\rm th}$ body,  
`$|\cdot|$' denotes the euclidean norm and `$\dot{\phantom q}$' denotes time derivative. 

\nl
In \cite[Chapter III, p. 125]{arnold}, V.I. Arnold made the following\footnote{The integer $n$ in Arnold's statement corresponds to the above $(1+n)$.}:

\nl
{\bf Arnold's statement:} ``{\sl In the $n$--body problem there exists a set of initial conditions having a positive Lebesgue measure and such that, if the initial positions and velocities belong to this set, the distances of the bodies from each other will remain perpetually bounded.}''

\nl
As well known, such statement solves a fundamental problem considered, for several centuries,  by astronomers and mathematicians. However, Arnold considered in details only the planar three--body case\footnote{Few lines after the above reported statement in \cite[Chapter III, p. 125]{arnold}, Arnold says: ``{\sl We shall consider only the plane three--body problem in detail. {\rm [$\cdots$]} In the final section a brief indication is given of the way in which the fundamental theorem of Chapter IV is applied in  the investigation of the planetary motions in the plane and spatial many--body problems.}''.} and it appears that his indication for extending the result to the general case contains a flaw; compare end of \S 1.2, p. 1524 in \cite{fejoz}.

\nl
A complete general proof of Arnold's statement was given only in 2004, when J.~F\'ejoz, completing the work of M. Herman,  proved the following\footnote{For a more detailed statement, see, footnote~\ref{thm.F} below.} 

\begin{teo}[\bf{Arnold, Herman, F\'ejoz \cite[\S 1.2, p. 1523, T{\scriptsize H\'EOR\`EME} 1]{fejoz}}]\label{AHF}$\phantom{AAAAAA}$
Si le maximum $\e = \max\{m_j/m_0\}_{j=1,...,n}$ des masses des plan\`etes rapport\'ees \`a la masse du soleil est suffisamment petit, les \'equations {\rm (1)} admettent, dans l'espace des phases au voisinage des mouvements k\'epl\'eriens circulaires et coplanaires, un ensemble de mesure de Lebesgue strictement positive de conditions initiales conduisant \`a des mouvements quasip\'eriodiques.
\end{teo}

\noindent
The beautiful proof of this result given in \cite{fejoz} (see also \cite{fejoz.w}) relies, 
on one side, on the elegant  $C^\infty$ KAM theory worked out by Herman (\S~2$\div$5 in \cite{fejoz}), and, on the other side, on the analytical celestial mechanics worked out, especially,  by Poincar\'e and clarified and further investigated in Paris in the late 1980's by A. Chenciner and J.~Laskar in the {\sl Bureau des Longitudes}\footnote{Compare, e.g., the {\sl Notes Scientifiques et techniques du Bureau des Longitudes} {\bf S 026} and {\bf S 028} by, respectively, Chenciner \& Laskar and Chenciner. 
}
and later by Herman himself.
  
\nl
The invariant tori associated to the motions provided by Theorem~\ref{AHF}, in view  of the just mentioned KAM tools, are $C^\infty$. Now, since the many--body problem is formulated in terms of real--analytic functions, it appears somewhat more 
natural to seek for {\sl real--analytic  invariant manifolds}. This is the problem addressed in this paper. In particular,  we shall give a new proof of  Arnold's statement, proving the following

\begin{teo}\label{thm.an}
If  $\e = \max\{m_j/m_0\}_{j=1,...,n}$ is small enough, there exists  a strictly positive measure set of initial conditions for the $(1+n)$--planetary problem {\rm (1)}, whose time evolutions lie on real--analytic  Lagrangian tori in the   $6n$--dimensional phase space  $$\M:=\Big\{(q,p)\in\R^{6(1+n)}: \ q_j\neq q_k\, ,\, \forall\, j\neq k\ {\rm and}\ \sum_{j=0}^n  p_j=0=\sum_{j=0}^n m_j q_j\Big\}\ ,$$ 
endowed with  the restriction of the standard symplectic form $\sum_{j=0}^n dq_j\wedge dp_j= \sum_{0\le j\le n\atop 1\le k\le 3} dq_{j,k}\wedge dp_{j,k}$.
\end{teo} 

\begin{oss}\label{oss1}\rm 
{\small
Let us collect here a few observations concerning the above statements and respective proofs.
\begin{itemize}
\item[(i)] The proof of Theorem~\ref{thm.an} given below is similar in strategy to that in \cite{fejoz} but technically different and it is based on an {\sl analytic} (rather than smooth) KAM theory for properly--degenerate Hamiltonian systems (see, also, point (iv) below).
On the other hand, 
it is conceivable -- notwithstanding  the presence of strong degeneracies  (see point (ii) and (iii) below) -- to prove regularity and uniqueness results for the planetary problem  so as to  deduce that the invariant tori in \cite{fejoz} are indeed analytic, a fact which does not follow from our proofs\footnote{Recent interesting progresses in the study of uniqueness of invariant Lagrangian manifolds appeared in \cite{BT} and, especially, in \cite{FGS}; however, as far as regularity is concerned, to the best of our knowledge, the only complete proven statement is Theorem 4, \S 4, p. 34 in \cite{salamon}, which  covers only the $C^\infty$ nondegenerate case.}. 

\item[(ii)] The evolution equations (\ref{Newton})  are Hamiltonian and admit seven integrals, namely, the Hamiltonian (energy)  $ H:=\sum_{j=0}^n \frac{|p_j|^2}{2m_j}- \sum_{0\le k<j\le n} \frac{m_jm_k}{|q_k-q_j|}$, the three components of the total linear momentum 
$ M:=\sum_{j=0}^n p_j$ and the three components of  the total angular momentum $C:= \sum_{j=0}^n p_j \times q_j$, where `$\times$' denotes the usual skew vector product in $\R^3$.  As a reflection of the invariance of Newton's equation (\ref{Newton}) under changes of inertial reference frames, the Hamiltonian system associated to the $(1+n)$--body problem may be studied on the symplectic, invariant $6n$--dimensional manifold $\M$ defined above,
where, besides the total linear momentum, also the coordinates of the barycenter of the system vanish (``reduction of the total linear momentum'').
However, the reduced $(1+n)$--body Hamiltonian still admits, besides the energy,  three integrals given by the components of $C=(C_x,C_y,C_z)$. Incidentally, such integrals are not commuting since, if $\{\cdot,\cdot\}$ denotes the natural Poisson bracket on $\M$, one has the cyclical relations $\{C_x,C_y\}=C_z$, $\{C_y,C_z\}=C_x$ and $\{C_z,C_x\}=C_y$; but for example $|C|^2$ and $C_z$ are two independent, commuting integrals. 

\item[(iii)] 
The reasons why, notwithstanding the development of KAM theory in the  early 1960's,  it took so long to give a complete proof of Arnold's statement are technical in nature and are related to the strong degeneracies of the planetary problem (degeneracies, which are related  to the abundance of integrals mentioned in the preceding item). 
The planetary $(1+n)$--body problem is  perturbative, the  unperturbed limit being
obtained by considering $n$ decoupled two--body problems formed by the Sun and the $j^{\rm th}$ planet. Now, the two--body problem in space is a three--degrees--of--freedom problem, but, once it is put into (Delaunay) action--angle  variables, it depends only on one  action  (the action $L$ proportional to the square root of the semi--major axis of the   Keplerian ellipse on which the two bodies revolve).
Systems of this kind are called {\sl properly degenerate} and standard KAM theory does not apply. This  difficulty, however,  was  overcome by Arnold -- essentially by refined normal form theory -- in the case of the planar three--body case, to which he could apply his ``fundamental theorem'' \cite[chapter IV]{arnold}. Indeed, Arnold's approach, in view of Jacobi's reduction of the nodes,  could   be extended  \cite{robutel} to the spatial  
three--body case ($n=2$)   but  not to the general case (spatial, $n>2$).
Indeed, when $n>2$, in space, a nice symplectic reduction (corresponding to the reduction of the nodes for $n=2$) is not yet known.

\noindent
Furthermore (but not independently), in higher dimension, there appear  two {\sl secular resonances} (see  Eq. (\ref{resonances}) below), which prevent direct application of any kind of KAM machinery. We mention that the way 
we overcome, here,  this last difficulty is slightly different from that  used in \cite{fejoz}: 
roughly speaking,  in \cite{fejoz} it is  introduced a modified Hamiltonian, which is then  considered  on  the symplectic submanifold   of vertical total angular momentum; here, we consider, instead,  an extended phase space by adding an extra degree--of--freedom and consider on it a modified nondegenerate Hamiltonian. 

\item[(iv)]  The main technical tool for us is the analytic KAM theory for weakly nondegenerate systems worked out by \Russ \ in \cite{russmann}; the main results of \Russ's theory (in the case of Lagrangian tori)  are recalled in \S~\ref{sub:R} (see, also, Lemma~\ref{lemmae_0} in \S~\ref{subsec:conc}).  The extension of this theory to properly degenerate systems is explained in \S~\ref{ssezteodegenere} and proved in \S~\ref{sezapplteoruss} (which constitutes the longest and most technical part of the paper). 
Finally, in \S~\ref{sezfinale}, using several results reported in \cite{fejoz}, the proof of Theorem~\ref{thm.an} is given. 

\item[(v)]Finally, we mention very briefly a few   questions related to the context considered here:

\begin{itemize}
\item Describe, in detail, the motions that take place on the Lagrangian tori. Let us clarify this point. 
From the proof of Theorem~\ref{thm.an} given below,  in view of the indirect argument used (\cite[Lemma 82, p. 1578]{fejoz}) we cannot conclude that the `true' motion is quasi--periodic; 
on the other hand, using  different arguments, Arnold  and F\'ejoz  
say that the motion, in the general case, is quasi--periodic and takes place on $(3n-1)$--dimensional tori\footnote{\label{thm.F} Arnold (\cite[p. 127]{arnold}): ``{\sl Thus, the Lagrangian motion is conditionally periodic and to the $n_o$ `rapid' frequencies of the Keplerian motion are added $n_o$ (in the plane problem) or $2n_o-1$ (in the space problem) `slow' frequencies of the secular motions}''.\\
F\'ejoz (\cite[p. 1566]{fejoz}):  {\bf T{\scriptsize H\'EOR\`EME} 60.} 
{\sl Pour toute valeur des masses $m_0$, $m_1$,..., $m_n >0$ et des demi grands
axes $a_1 > \cdots > a_n > 0$, il existe un r\'eel  $\e_0 > 0$ tel que, pour tout  $\e$  tel que
$0 < \e  < \e_0$, le flot de l'hamiltonien $F$ (d\'efini en (28)) poss\`ede un ensemble de mesure
de Lebesgue strictement positive de tores invariants de dimension $3n-1$, de classe $C^\infty$,
quasip\'eriodiques et  $\e$-proches en topologie $C^0$ des tores k\'epl\'eriens de demi grands axes
$(a_1, . . . ,a_n)$ et d'excentricit\'es et d'inclinaisons relatives nulles; de plus, quand   tend vers
z\'ero la densit\'e des tores invariants au voisinage de ces tores k\'epl\'eriens tend vers un.}}. Moreover,  in the spatial three--body case ($n=2$) the Lagrangian tori are actually $4$--dimensional (not $5=3n-1$) and the number of {\sl independent} frequencies is 4 (compare \cite{robutel}). 

\item Find a `good' set of analytic symplectic variables for the general spatial many--body problem.

\item Give asymptotic (as $\e\to 0$) estimates on the measure of Lagrangian invariant tori.

\item Apply some of the above result to a subsystem of the Solar system (for some progress in this direction, see \cite{CC}). 

\end{itemize}

\end{itemize}

}
\end{oss}

\section{Analytic Lagrangian tori for properly degenerate systems}
\label{russmanndegenere}
In this section we first recall a  result  due to \Russ\ concerning analytic perturbations of weakly nondegenerate Hamiltonian systems (\S~\ref{sub:R}) and then show how such result may be used to give an analytic version of Herman's $C^\infty$ KAM theorem on properly degenerate systems (i.e., nearly--integrable systems, which when the perturbation parameter vanishes depend on less action variables than the number of degrees of freedom). The statement of the analytic theorem for properly degenerate systems is given in \S~\ref{ssezteodegenere} and its proof in \S~\ref{sezapplteoruss}.

\subsection{\Russ's theorem for weakly  nondegenerate systems}
\label{sub:R}
We start with fixing some notation. 

\begin{itemize}{\footnotesize
\item If  $a,b \in \R^n$ then $\langle a, b \rangle := \sum_{i=1}^n a_i b_i$ 
and $|a| := {|a|}_2 := { \langle a, a \rangle }^\frac{1}{2}$; 

\item  if $g$ is a $\mu$--times continuously differentiable function ($\mu\in\N$) from an open set $B \subset \R^n$ to $\R^m$,
the $\mu$--th (tensor) derivative of $g$ in $b \in B$ is denoted by 
$(a_1, \dots a_\mu) \rightarrow \partial^\mu g(b) (a_1, \dots, a_\mu)$, $a_j \in \R^n, j=1, \dots \mu$; if $a_1=\cdots=a_\mu$,  we shall write $\partial^\mu g(b) (a)^\mu$

\item 
$|\partial^\mu g(b)|  :=  \max_{a \in \R^n, |a| = 1 } | \partial^\mu g(b) (a, \dots, a) |$ and 
${|\partial^\mu g|}_A := \sup_{b \in A} |\partial^\mu g(b)|$;

\item 
$C^\mu (B, \R^m)$ will denote  the Banach space of all $\mu$--times continuously differentiable functions
$g: B \rightarrow \R^m$ with bounded derivatives up to order $\mu$, endowed with the norm ${|g|}_B^\mu = \sup_{0 \leq \nu \leq \mu} {|\partial^\nu g|}_B  <  \infty$.
}
\end{itemize}

\nl
The  key   notion of nondegeneracy is the following. 

\begin{defin}[\Russ \, nondegeneracy condition\footnote{This terminology seems to be, nowadays, standard (see, e.g., \cite{sevryuk}); however many authors, besides \Russ, contributed to its formulation; among them: Arnold, Margulis, Pyartli, Parasyuk, Bakhtin, Sprindzhuk and others.}]                                                       
\label{defnondegR}
A real--analytic function 
$$\o : y\in B \subset \R^n \longrightarrow \o (y) = (\o_1(y) , \dots , \o_m(y)) \in \R^m$$
is  called {\bf R--nondegenerate} if 
$B$ is a non--empty open connected set in $\R^n$ and if
for any  $ c = (c_1, \dots , c_m)  \in  \R^m  \smallsetminus \{ 0 \} $ 
one has
\begin{equation*}
y \longrightarrow \langle c, \o \rangle := \sum_{i=1}^m c_i \o_i \neq 0
\end{equation*}
or equivalently  if the range $\o(B)$ of $\o$  does not lie in any $(m-1)$--dimensional linear subspace of $\R^m$.
We call $\o$ {\bf R--degenerate} if it is not R--nondegenerate.
\end{defin}

\nl
The following lemma is a simple consequence of R--nondegeneracy and analyticity:

\medskip
\begin{lem}                                                                     \label{lemmaindicenondeg}
Let $\o : B \subset \R^n \longrightarrow \R^m$ be R--nondegenerate.
Then for any non empty compact set $ \K \subset B $  there exist numbers
$\mu_0 = \mu_0 (\o , \K) \in \Z_+  $ and $\beta = \beta (\o , \K) > 0$  such that
\begin{equation}
\label{indiceamount}
\max_{0 \leq \mu \leq \mu_0} \left| \partial^\mu_y { \langle c, \o(y) \rangle }^2 \right|  \geq \beta  \,\,\,  ,
							\for c \in \mathcal{S}^{m-1} , \for y \in \K 
\end{equation}
where $\mathcal{S}^{m-1} := \{ c \in \R^m \, : \, {|c|}_2 = 1 \}$.
\end{lem}

\medskip

\nl
For the proof see Lemma~18.2 on page 185 of \cite{russmann}.

\nl
In view of Lemma~\ref{lemmaindicenondeg} one can give the following 

\begin{defin}
\label{defindexamount}
Let $\K$ and $B$ as in the preceding lemma and
let  $\o : y \in B \longrightarrow \R^m$ be a real--analytic and R--nondegenerate function. 
We define $\mu_0 (\o , \K) \in \Z_+$,  the {\bf index of nondegeneracy} of $\o$ with respect to $\K$,
as the smallest positive integer such that
\begin{equation}
\label{beta}
\beta := \min_{y \in \K, \, c \in \mS^{m-1} } \, \max_{0 \leq \mu \leq \mu_0} 
														\left| \partial^\mu { \langle c, \o(y) \rangle }^2 \right|  > 0 \, .
\end{equation}
The number $\beta = \beta (\o, \K)$ is called {\bf amount of nondegeneracy} of $\o$ with respect to $\K$.
\end{defin}
\medskip

\begin{oss}\rm
\label{ossRnondeg}
If a real--analytic function $\o : B \rightarrow \R^m$ admits the existence of $\mu_0$
and $\beta$ as in (\ref{beta}), for some compact $\K$ containing an open ball,
then it is R--nondegenerate.
\end{oss}\rm 

\nl
The following result, which concerns the existence of maximal (Lagrangian) tori only, 
is a particular case of the main Theorem in \cite{russmann},
where also lower dimensional tori are treated\footnote{See Theorem $1.7$ in \cite[p. 127]{russmann}.
We refer to \cite{miatesi} for more details 
on how to obtain Theorem \ref{teorussmassimale} from the general results in \cite{russmann}.}.

\begin{teo}[\bf{\Russ, 2001}]                         \label{teorussmassimale}
Let $\mathcal{Y}$ be an open connected set of $\R^n$ and $\T^n$ the usual $n$--dimensional torus
$\R^n/ 2 \pi \Z^n$. 
Consider a real--analytic Hamiltonian
\begin{equation*}
H (x , y )  =  h ( y )  +  P (x, y)
\end{equation*}
defined for $ (x, y)  \in \T^n \times \mathcal{Y} $ endowed with the standard symplectic form $dx \wedge dy$.
Let $\K$ be any compact subset of $\mathcal{Y}$ with positive n--dimensional Lebesgue measure $\meas_n\, \K >0$ 
and fix $0 < \e^\star < \meas_n\, \K$.
Let $\A$ be an open set in $\C^n / 2 \pi \Z^n \times \C^n$ on which $H$ can be analytically 
extended and such that $ \T^n \times \K \subset \A$.
Assume that the frequency application $\o := \nabla h$ is R--nondegenerate on $\mathcal{Y}$;
let $\mu$ be any integer greater or equal to $\mu_0 (\o, \K)$ 
(the index of nondegeneracy of $\o$ with respect to $\K$) and let $\beta$ be as in (\ref{beta})
with $\mu_0$ replaced by $\mu$.

\nl
Then, for any fixed $\tau > n \mu$, there exist 
$\e_0 = \e_0 ( \e^\star, n, \mu, \beta, \tau, \o, \K) > 0$ 
and 
$\g = \g (\e^\star, n, \tau, \mu, \beta, \o, \K) > 0$
such that if
\begin{equation}
\label{tagliaperturbazione}
{|P|}_\A := \sup_{\A} |P| \leq \e_0
\end{equation}
the following is true.
There exist a compact set 
\begin{equation}
\label{misuraK}
\K^\star \subset \K \quad{\rm  with}\quad \meas_n\, \K^\star > \meas_n\, \K - \e^\star
\end{equation}
and a Lipschitz mapping
\begin{equation*} 
X :   ( b , \xi , \eta )  \in  \K^\star  \times \T^n  \times  \mathcal{U}  \longrightarrow  \T^n \times \Y   \, ,
\end{equation*}
where   $ \mathcal{U} $  is an open neighborhood of the origin in $\R^n$, such that:
\begin{itemize}
      					\item[(i)] the mapping
                \begin{equation*}
                (\xi ,  \eta)   \longmapsto  (x, y)  =  X (b, \xi, \eta)
                \end{equation*}
                defines, for every $b \in \K^\star $,  a real--analytic symplectic transformation\footnote{
                I.e. it preserves the symplectic form $dx \wedge dy$.} close to the identity on
                $ \T^n \times \mathcal{U}$;
                
                \item[(ii)] the map
                \begin{equation*}
                (b, \xi) \in \K^\star \times \T^n\longrightarrow X_0 (b, \xi) := X (b, \xi, 0)
                \end{equation*}
                is a bi--Lipschitz homeomorphism;

                \item[(iii)] the transformed Hamiltonian $H^\star := H \circ X$ is in the form\footnote{
                Here and in what follows
                $f(\eta) = O( g(\eta) ) $ means that there exists a constant $C$ such that
                $|f (\eta)| \leq C |g(\eta)| $ for small enough $\eta$.}:
                \begin{equation*}
                 H^\star (b, \xi, \eta)  =  
                 h^\star ( b )  + \langle \o^\star (b) ,  \eta \rangle  +  O ( { | \eta | }^2  )
                \end{equation*}
                for every $b \in \K^\star$ and $ (\xi , \eta) \in \T^n \times \mathcal{U}$;

                \item[(iv)] the new frequency vector $\o^\star$ satisfies for all $ b $  in $ \K^\star $  the Diophantine inequality
                \begin{equation}
                \label{disdiofantine}
                |  \langle  k , \o^\star (b)   \rangle  | \geq  \frac{\g}{ {| k |}_2^\tau }     ,
                                                        \for  k \in \Z^n \smallsetminus \{  0 \} \, .
                \end{equation}
\end{itemize}
\end{teo}
\medskip

\begin{oss}\rm 
\begin{itemize}\item[(i)] From Theorem \ref{teorussmassimale} we immediately obtain that for any $b \in \K^\star$
			the $n$--dimensional tori
			\begin{equation*}
			\mathcal{T}_b := X_0 (b, \T^n)
			\end{equation*}
			are invariant for $H$ and the $H$--dynamics
			is analytically conjugate to $\xi \rightarrow \xi + \o^\star(b) t$.
			Furthermore, as it follows from (\ref{misuraK}) and point {\sl (ii)},  the measure of $ \, \cup_{b \in \K^\star} \mathcal{T}_b $ 
			is proportional to $(\meas_n\, \K - \e^\star){(2 \pi)}^n$ and, hence, tends to the full measure linearly when $\e^\star$ tends to 0.

\item[(ii)] A (technical) difference between \Russ's Theorem and the formulation given above in Theorem \ref{teorussmassimale} is
			the choice of ${\mu}$ as any integer greater or equal than the actual index of nondegeneracy of $\o$,  while in \cite{russmann} $ \mu$ is chosen equal to the index of nondegeneracy of $\o$. In fact, it is easy to check\footnote{See, e.g., \cite{miatesi}.} that \Russ's theorem holds in this slightly more general case, which will however be important in our applications.

Another difference of Theorem \ref{teorussmassimale}  above with respect to \Russ's original formulation, concerns the way the small divisors are controlled. 
\Russ \, uses a very general approach based upon  ``approximation functions'';
however, such approach is too general for our application and cannot be applied directly.  Nevertheless, it is easy to follow a more classical approach\footnote{Compare, again, \cite{miatesi}.} based upon Diophantine inequalities of the form (\ref{disdiofantine}), which will be good enough for the application to properly degenerate systems; compare  also remark \ref{ossPhi}, (ii)  below.
\end{itemize}
\end{oss}\rm

\subsection{A KAM Theorem for  properly--degenerate systems }
\label{ssezteodegenere}
Let $d$ and $p$ be positive integers; let $\mathcal{B}$ an open set in $\R^d$,  $\mathcal{U}$ some open neighborhood of the origin in $\R^{2p}$
and $\e$ a ``small'' real parameter. Consider a Hamiltonian function $H_\e$ of the form
\begin{equation}																
\label{Hamiltonianainiziale}
H_\e (\varphi,I,u,v) = h(I) + \e f (\varphi,I,u,v) \, ,
\end{equation}
real--analytic for
\begin{equation*}
(\varphi, I , (u, v) ) \in \T^d \times \mathcal{B} \times \mathcal{U} =: \mathcal{M}
\end{equation*}
where $\mathcal{M}$ is endowed with the standard symplectic form
\begin{equation*}
d \varphi \wedge  d I  +  du \wedge dv \, .                                                     
\end{equation*}
The ``perturbation'' $f$ is assumed to have the form 
\begin{equation}
\label{formaf}
\left\{
\begin{array}{l}
f (\varphi, I,u, v ) = f_0 (I,u,v) + f_1 (\varphi, I,u,v)  
\,\, , \,\,\,\, \int_{\T^d} f_1 (\varphi,I,u,v) \, d \varphi = 0
\\
\\ 
\displaystyle f_0 (I,u,v) = f_{00}(I) + \sum_{j=1}^p \Omega_j (I) \frac{u_j^2 + v_j^2}{2}
                                                        +  O \left( {|(u,v)|}^3 \right)  \, .
\end{array}
\right.
\end{equation}

\nl
Observe that the Hamiltonian $h + \e f_0$  possesses for every $\bar{I} \in \mathcal{B}$ the invariant isotropic 
(non--Lagrangian) torus
\begin{equation*}
\mathcal{T}^d_{\bar{I}} := \T^d  \times \{ \bar{I} \}  \times \{ 0 \}  \subset \mathcal{M} 
\end{equation*}
with corresponding quasi--periodic flow
\begin{equation*}
\varphi (t) = \Big( \partial_I h (\bar{I}) +
                        \e \partial_I f_{00} (\bar{I}) \Big)t +  \varphi_0
                        \qquad  I(t) \equiv \bar{I}  \qquad (u(t), v(t)) \equiv 0      \,  .
\end{equation*}
The purpose is to find Lagrangian invariant tori for $H_\e$ close to $(d+p)$--tori of the form
\begin{equation}      
\label{d+ptori}
\mathcal{T}^{d + p}_{\bar{I} , w} = \T^d \times \{ \bar{I} \}  \times
                        \{ (u,v) \in \R^{2p} \, , \, {|(u_j ,v_j)|}^2 = 2w_j  \, ,\for j=1 , \dots , p\}
\end{equation}
for $\bar{I}$ in $\mathcal{B}$ and 
$w \in { \left( \R_+ \right) }^p$ small.

\begin{teo}
\label{teoremadegenere}
Consider a real--analytic Hamiltonian function $H_\e$ as in (\ref{Hamiltonianainiziale}) 
and (\ref{formaf}), and assume that the ``frequency map''
\begin{equation}                                                                                
\label{applicazionefreq}
I \in \mathcal{B}  \longrightarrow  ( \o(I) , \Omega(I) ) := 
					( \nabla h(I), \Omega_1(I), \dots, \Omega_p (I) ) \in \R^d \times \R^p
\end{equation}
is R--nondegenerate.
Then, if $\e$ is sufficiently small, 
there exists a positive measure set of phase space points belonging to real--analytical, 
Lagrangian, $H_\e$--invariant tori, which are close to $ \mathcal{T}^{d + p}_{ \bar{I} , w}$ as in (\ref{d+ptori}) 
with $w_j = O( \e )$; furthermore, the $H_\e$--flow on such tori is  quasi--periodic
with Diophantine frequencies. 
%of the form $(\o + O(\e) , \e \Omega + O(\e^2))$. 
\end{teo}

\begin{oss}\label{oss*}\rm 
\begin{enumerate}
	\item[(i)] The above Theorem may be viewed as the real--analytic version for Lagrangian  tori of the $C^\infty$   KAM Theorem by M. Herman contained in \cite{fejoz} 
	(see, in particular, Theorem 57, page 1559). Under stronger nondegeneracy assumptions the above theorem corresponds to the ``Fundamental Theorem'' in \cite{arnold}.
			
	\item[(ii)] The word ``properly--degenerate'' refers to the fact that for $\e = 0$ the 
	Hamiltonian $H_0$ depends on $d$ action variables, while the number of degrees of freedom is $d+p>d$. In particular, the 
	tori constructed in Theorem~\ref{teoremadegenere}, as $\e\to 0$, 	degenerate into lower dimensional (non Lagrangian) tori $\mathcal{T}^d_{\bar I}$.
	
	\item[(iii)] The natural symplectic variables for the KAM theory of the Hamiltonian $H_\epsilon$ are $(I, \varphi)$ and (rather than the cartesian variables $(q,p)$) the symplectic action--angle variables $(w, \zeta)$,  where $w_j = \frac{u_j^2 + v_j^2}{2}$ for $j=1, \dots, p$ and $\zeta_j$ is the angle of the circle\footnote{Compare Eq. (\ref{coordinatepolari}) below, where $w$ is related to $\rho$ by $w = \rho^0 + \rho$.} $|w_j|=$ const.
Indeed, Theorem~\ref{teoremadegenere} has, in terms of such variables, a natural reformulation, which gives a deeper insight into the structure of the invariant tori\footnote{In  reformulating Theorem~\ref{teoremadegenere} in terms of the variables $(\zeta,w)$ we shall often use the same symbols used above. The Proof of Theorem~\ref{teo5} will not be explicitly given since it follows easily from the proof of Theorem~\ref{teoremadegenere}.}:

\begin{teo}\label{teo5}
Let $H_\e (\varphi,I,\zeta,w) = h(I) + \e f (\varphi,I,\zeta,w)$ be real--analytic for
$(\varphi, I , \zeta, w ) \in \T^d \times \mathcal{B} \times \T^p\times \{w\in\R^p: 0<|w_j|<r\} =: \mathcal{M}$ for some open set $\mathcal{B}\subset \R^d$ and $0<r$; $\mathcal{M}$ is endowed with the symplectic form  $d \varphi \wedge  d I  +  d\zeta \wedge dw$. The perturbation $f$ is of the form $f=f_0(I,\zeta,w)+f_1$ with $f_1$ having vanishing $\varphi$--mean value over $\T^d$; furthermore  $f_0$ has the form
$f_0= f_{00}(I) + \langle \Omega (I), w\rangle+  o(|w|)$. Then, if the frequency map
$I\in \mathcal{B}\to (\omega,\Omega):=(\partial_Ih(I),\Omega(I))\in\R^d\times\R^p$ is R--nondegenerate, and $\e$ is small enough, there exists a positive measure set of phase space points belonging to real--analytical, 
Lagrangian, $H_\e$--invariant tori, which have the following parametrization:
		$$
		\left\{
		\begin{array}{l}
		\varphi  =  \theta + \tilde{\varphi} (\theta, \psi)
		\\
		\\
		I  =  \bar{I} + \tilde{I} (\theta, \psi)
		\\
		\\
		\zeta  =  \psi + \tilde{\zeta} (\theta, \psi)
		
		\\
		\\
		w  =  \bar{w}  +  \tilde{w} (\theta, \psi)
		\end{array}
		\right.
		$$
		where $\bar{w}$ is a constant vector of norm $2 \e$ and
		$\tilde{\varphi}, \tilde{I}, \tilde{\zeta}$ and $\tilde{w}$ real--analytic functions 
		for $(\theta, \psi) \in \T^d \times \T^p$ (with range, respectively, in $\T^d, \R^d, \T^p$ and $\R^p$)
		with
		$$
		\left\{
		\begin{array}{l}
		\tilde{I}  =  O \left( \e  {\left( \log \e^{-1}  \right)}^{- (\tau_0 + 1) } \right)
		\\
		\\
		\tilde{w}  =  O \left( \e^\frac{\nu+1}{2} \right)
		\\
		\\
		\tilde{\varphi} \,\,  ,  \, \tilde{\zeta} = O( \e )  
		\end{array}
		\right.
		$$
		for suitable $\nu \geq 4$ and\footnote{
		See equations (\ref{Hconiugata}) and (\ref{sa_1tau_0}) below.} $\tau_0 \geq d+p$.
		Moreover, if $(\o, \Omega)$ is the above frequency map,  the $H_\e$--flow on such invariant  tori is conjugated to
		\begin{equation*}
		(\theta, \psi)  \longrightarrow \left( \theta + \tilde{\o} t ,  \psi + \e \tilde{\Omega} t \right)
		\end{equation*}
		for a suitable Diophantine vector $(\tilde{\o}, \tilde{\Omega})$ satisfying
		\begin{equation*}
		| \tilde{\o} - \o | \,\,  ,  \,\,\,  | \tilde{\Omega} - \Omega |  =   O  ( \e )   \,\,  .
		\end{equation*}
	\end{teo}
\end{enumerate}	
\end{oss}

\subsection{Proof of Theorem \ref{teoremadegenere}}
\label{sezapplteoruss}
First of all, let us introduce some  notation and make quantitative the assumptions of Theorem~\ref{teoremadegenere}.
\begin{itemize}
	\item For $\d>0$, $d\in\N$, $A \subset \R^d$ or $\C^d$   we denote
	\begin{eqnarray}
	\label{notazioneinsiemi0}
	B^d (x_0, \d) & :=& \{ x \in \R^d \, : \, |x - x_0| < \d\}\ ,\qquad (x_0\in\R^d)\ ,
\\ \  \nonumber	
	D^d (x_0, \delta)&:=&\{ x  \in \C^d \, : \, |x  - x_0| < \delta  \}  \ ,\qquad (x_0\in\C^d)\ ,
\\ \nonumber
	\T^d_\delta & := &  \{ x \in \C^d \, : \,  | \Im  x_j | < \delta, \Re x_j \in \T  \, , \for j = 1 \dots d  \}
	\\ \nonumber
	A + \delta & := & \bigcup_{ x \in A }  D^d (x, \delta) 
	\label{notazioneinsiemi}
	\end{eqnarray}

	\item We may assume that $H_\e$ in (\ref{Hamiltonianainiziale}) and (\ref{formaf})
	can be holomorphically extended for 
	\begin{equation}
	\label{dominiooloH} 
	(\varphi, I , (u, v) ) \in \T^d_\s \times (\mathcal{B} + r_0) \times  (\mathcal{U} + r_1) =: \M_\star  \, .
	\end{equation}
	In particular $H_\e$ is real--analytic on $\T^d \times B^d (I_0, s) \times B^{2p} (0 ,r_1)$
	for any   $I_0$ in $\mathcal{B}$ and $s < r_0$.
	Moreover, we shall denote
	\begin{equation}
	\label{M_0}
	M_0 :=  \sum_{k \in \Z^d} \left(\sup_{ (\mathcal{B} + r_0) \times (\mathcal{U} + r_1) } 
				|f_k(I,u,v )| \right) e^{ {|k|}_1 \s}
	\end{equation}
	as the ``sup--Fourier'' norm of $f$ and let
	\begin{equation}
	\label{normaoO}
	M_1  := \sup_{ I \in \B + r_0 } | ( \o (I), \Omega(I) ) |    \, .
	\end{equation}
\end{itemize}

\nl
The proof of Theorem \ref{teoremadegenere} is   based upon {\bf two preliminary steps}: \begin{itemize}
\item[\bf  1] computation of a suitable normal form for $H_\e$;
\item[\bf  2]  quantitative estimates on the amount of the nondegeneracy of the normal form.
\end{itemize}

\subsubsection{Step 1: Normal forms for properly--degenerate systems}
\label{sezdimlemma}

\begin{pro}                                                                                      
\label{teoremaconiugazione}
Fix  an integer $\nu \geq 4$. Then, there exists $m>d$ (depending on $\o$), 
and, for  $\e$  small enough, a point $I_0\in \mathcal{B}$ and 
a real--analytic canonical transformation\footnote{
Symplectic up to rescalings.} $\Phi_\e $ such that the following holds.
Let 
\begin{equation}
\label{ordines}
s := O\Big( {\left( \log \e^{-1} \right)}^{ - m }  \Big)              
\end{equation}
then $B^d (I_0, s)\subset \mathcal{B}$ and 
$\Phi_\e:  (\vartheta, r, \zeta, \rho)\longrightarrow (\varphi, I, u, v)$ satisfies
$$ 
\Phi_\e : \T^d  \times  B^d (0 , s/5 )  \times \T^p  \times  B^p (0,\e)
          \longrightarrow  
          \T^d \times B^d (I_0, s) \times \mathcal{U}															\label{dominiconiugazione}
$$
and  $\hat{H}_\e := H_\e \circ \Phi_\e $ takes the form
\begin{eqnarray} 
\label{Hconiugata}
\hat{H}_\e  (\vartheta, r, \zeta, \rho ) &=& N_\e(r,\rho; \rho^0)+  \e^\nu  P_\e  (\vartheta, r, \zeta, \rho ; \rho^0)
\end{eqnarray}
with
\begin{equation}\label{Neps}
N_\e
:=
\frac{1}{\e} h (I_0 + \e r) +  \hat{g} (I_0 + \e r)   +     
\frac{1}{2}  \hat{\Omega} (I_0 + \e r)   \cdot  (\rho^0 + \e \rho) +    
Q_{\e , I_0 + \e r} (\rho^0 + \e \rho) 
\end{equation}
and:  $\rho^0$ in $ { \left(  \R_+ \right) }^p$ is some point having euclidean norm $2 \e$;
$ Q_{\e , I_0 + \e r} $ is a polynomial of degree $\nu-1$ starting with cubic terms;
$\hat{g}$ , $\hat{\Omega}$ and $P_\e$ are real--analytic functions.
Furthermore, one has
\begin{equation*}                                                                                %\label{hatO-O=e}
\sup_{ r \in B^d (0, s/5) } | \hat{\Omega} (I_0 + \e r) -  \Omega (I_0 + \e r) |  =  
										O \left( \e {\left( \log \e^{-1} \right)}^{2m - 1} \right)   \, .
\end{equation*}
\end{pro}
\nl {\bf Proof of Proposition \ref{teoremaconiugazione}.}
We start by  recalling  a measure theoretical result due to Pyartli (see \cite{pyartli} or  \cite[Theorem 17.1]{russmann}):
\begin{lem}[{\bf Pyartli}]
\label{teopyartli}
Let $\K \subset \R^d$ be a compact set; let $\theta \in (0,1)$ 
define $\K_\theta := \cup_{y \in \K} B^d (y, \theta)$. 
Let $g : \K_\theta \rightarrow \R$ be a real--analytic function satisfying
\begin{equation*}
\min_{y \in \K } \, \max_{0 \leq \nu \leq \mu_0} 
\left| \partial^\nu g(y) \right|  \geq  \beta
\end{equation*}
for some $\beta > 0$.
Then there exists $C = C (\mu_0, \beta, d, \K, \theta)$ such that
\begin{equation*}
\meas_d\, \left\{ y \in \K  \, : \,  | g(y) | \leq t  \right\} \leq 
													C {|g|}_{\K_\theta}^{\mu_0 + 1} t^\frac{1}{\mu_0}
\end{equation*}
for any $0 \leq t \leq \frac{\beta}{2 \mu_0 + 2}$.
\end{lem}

\nl
Pyartli's Lemma implies: 

\begin{lem}
\label{corpyartli}
Let $\K$ be a compact set with positive $d$--dimensional Lebesgue measure and let let $0 < \varepsilon^\star < \meas_d\, \K$; 
let $\o : \K_\theta \rightarrow \R^d$ be R--nondegenerate and let $\mu_0$ and $\beta$ be its
index and amount of nondegeneracy with respect to $\K$.
Let us denote by $\mathcal{D}^d_{\g_0, \tau_0}$ 
the set of Diophantine vectors in $\R^d$ with Diophantine constants $\g_0, \tau_0$, i.e. the set
\begin{equation*}
\mathcal{D}^d_{\g_0, \tau_0} := \left\{ \o \in \R^d  \, : \,  | \langle \o, k \rangle | 
					\geq \frac{\g_0}{ {|k|}^{\tau_0} } \, , \for k \in \Z \smallsetminus  \{ 0 \}   \right\}  \, .
\end{equation*}
Then, if $\g_0 $ is sufficiently small and $\tau_0 \geq d \mu_0$ one has
\begin{equation}
\label{stimadiof}
\meas_d\,  \left( \K \cap \mathcal{D}^d_{\g_0, \tau_0} \right) \geq \meas_d\, \K - \varepsilon^\star \, .
\end{equation}
\end{lem} 

\nl 
{\bf Proof of Lemma~\ref{corpyartli}.}  First of all observe that for any $m \in \Z_+$, $a \in \R^d$, $k \in \Z^d \smallsetminus \{ 0 \}$
and $b \in \K_\theta$
\begin{equation*}
\left|  \partial^m \langle \o(b) , k {|k|}^{-1}  \rangle  (a^m)  \right|  =
			%\left|  \sum_{i=1}^d \\partial^\nu \o_j (b) (a^\nu)  \frac{k_j}{|k|}   \right|  \leq
			\left|  \langle \partial^m \o(b) (a^m) ,  k {|k|}^{-1}  \rangle \right|  \leq
			\left|  \partial^m \o(b) (a^m)  \right|  \, ;
\end{equation*}
taking the $\sup$ over $|a| = 1$, $b \in \K_\theta$ and $0 \leq m \leq \mu$ we get
\begin{equation*}
{ \left| \langle \o, k {|k|}^{-1}  \rangle \right| }^\mu_{\K_\theta} \leq  
						{ | \o | }^\mu_{\K_\theta} < \infty 
\end{equation*}
for any $\mu \in \Z_+$.
Now we use this last inequality and Theorem \ref{teopyartli}, assuming $\g_0 \leq \frac{\beta}{2 \mu_0 + 2}$,
to estimate
\begin{eqnarray*}
 \meas_d\, \left( \K \smallsetminus \mathcal{D}^d_{\g_0, \tau_0} \right)  &=&  
					\meas_d\,  \bigcup_{k \in \Z^d \smallsetminus \{ 0 \}}  
					\left\{ b \in \K  \, : \,  | \langle \o(b), k \rangle | < \frac{\g_0}{ {|k|}^{\tau_0} }  \, \right\}
					\\
					\\
					&   \leq &\!\!  \sum_{k \in \Z^d \smallsetminus \{ 0 \}}
					\meas_d \left\{ b \in \K  \, : \,  
					\left| \langle \o(b), \frac{k}{|k|} \rangle  \right| < \frac{\g_0}{ {|k|}^{\tau_0 + 1} }  \, \right\}
					\\
					\\
					& \leq&  C (\mu_0, \beta, d, \K, \theta)  { | \o | }^{\mu_0+1}_{\K_\theta}  \g_0^\frac{1}{\mu_0}
					\sum_{k \in \Z^d \smallsetminus \{ 0 \}} \frac{1}{ {|k|}^\frac{\tau_0+1}{\mu_0}  }   \, .
\end{eqnarray*}
Since $\tau_0 \geq d \mu_0 $ this last sum converges and one has
\begin{equation*}
\meas_d\, \left( \K \smallsetminus \mathcal{D}^d_{\g_0, \tau_0} \right)  \leq   \bar{C} \g_0^\frac{1}{\mu_0}
\end{equation*}
for a suitable $\bar{C} = \bar{C} (\mu_0, \beta, d, \K, \theta, \o, \tau_0)$. 
Choosing $\g_0 \leq { \left(  \bar{C}^{-1} \e^\star \right) }^{\mu_0}$  we obtain  estimate (\ref{stimadiof}). \hfill$\blacksquare$

\medskip

\nl
Now consider  the real--analytic Hamiltonian $H_\e$ in (\ref{Hamiltonianainiziale}) and (\ref{formaf}).
Let $\nu_1, \nu_2 \geq 4$ be two integers to be later determined and set
\begin{equation}
\label{definizioneK_1}
K_1 := \frac{6}{\sigma} (\nu_1 - 1) \log \frac{1}{\e M_0}
\end{equation}
where $M_0$ is defined by (\ref{M_0}).
Lemma \ref{corpyartli} and the R--nondegeneracy of $\o$ assure the existence of $I_0 \in \B$ such that 
$\o (I_0)$ belongs to $\mathcal{D}_{\g_0, \tau_0}^d$ (for suitable $\g_0$ and $\tau_0$).
Then, from Taylor's formula it follows that\footnote{
For $k = (k_1, k_2, \dots, k_d) \in \Z^d$ we denote ${|k|}_1 := \sum_{i=1}^d |k_i|$; recall also the definition of  complex balls $D^d$ in (\ref{notazioneinsiemi0}).}
\begin{equation}
\label{ouniformenonrisonanza}
|\o (I) \cdot k| \ge \alpha_1 >0 \, , \for k \in \Z^d ,\, 
						0 < {|k|}_1  \leq  K_1 , \,  \for I \in D^d (I_0 , s)
\end{equation}
with\footnote{
This means that we can take $m = \tau_0 + 1$ in (\ref{ordines}).} 
\begin{equation}
\label{sa_1tau_0}
s := O \left( {\left( \log \e^{-1} \right)}^{ -(\tau_0 + 1) }  \right)
\qquad \mbox{and} \qquad
\alpha_1 :=  O \left( {\left( \log \e^{-1} \right)}^{ - \tau_0 }  \right)  \, .
\end{equation} 
Furthermore, we can assume that there exists $\alpha_2$ (independent of $\e$) such that
\begin{equation}
\label{Ouniformenonrisonanza}
|\Omega (I) \cdot k| \geq \alpha_2 > 0  \, , \for k \in \Z^p ,\, 
						0 < {|k|}_1  \leq  \nu_2 , \,  \for I \in D^d (I_0, s) \, .
\end{equation}

\nl
Next, we want to average $H_\e$ over the ``fast angles'' $\varphi$ up to order $\nu_1$. To do this we shall apply the following classical ``averaging lemma'', whose proof 
can be found in \cite[Appendix A, p. 110]{BCV}.

\begin{lem}[\bf Averaging Lemma\footnote{
Lemma \ref{lemmamedia} can be immediately derived from Proposition A.1 in \cite{BCV} with the following correspondences:
$\alpha_1 = \alpha$ for $\alpha_1$ as in (\ref{ouniformenonrisonanza}), $K_1 = K$ for $K_1$ as in (\ref{definizioneK_1}), 
$\e M_0 = \varepsilon$ for $M_0$ as in (\ref{M_0}), 
$ s = r, d$ for $s$ as in (\ref{ordines}) and (\ref{ouniformenonrisonanza}),
$\{ 0 \} = \Lambda$ and $\e f (I,\varphi,u,v) $ in (\ref{Hamiltonianainiziale}) is just $f(u,\varphi)$ in \cite{BCV};
as a result one has that $\e f_0 + \tilde{g}$ and $\tilde{f}$ 
are respectively given by $g$ and $f_\star$ in \cite{BCV} with estimates (\ref{stimeteomedia}) holding
in view of the previous correspondences.
}]
\label{lemmamedia}
Let $H_\e, M_0, \sigma, \alpha_1$ and $s$ be as above.
Assume (\ref{ouniformenonrisonanza}) holds with $K_1$ as in (\ref{definizioneK_1}).
Then, if $\e$ is small enough,
there exists a real--analytic symplectic transformation $\Phi^1_\e : ( \tilde{\varphi}, \tilde{I}, \tilde{u},\tilde{v} ) \to  (\varphi, I, u, v) $ mapping 
$$\mathcal{M}_1 := 
\T^d_\frac{\s}{6} \times D^d (I_0 , \frac{s}2)  \times D^{2p} (0, \frac{r_1}2) 
\stackrel{\Phi^1_\e}{\to}  \mathcal{M}_0 := \T^d_\s  \times D^d (I_0 , s)  \times D^{2p} (0, r_1) 
$$
that casts $H_\e$ into the Hamiltonian
\begin{equation*}
H^1_\e := H_\e \circ \Phi^1_\e = h + \e f_0 + \tilde{g} + \tilde{f}
\end{equation*}
where  $ \tilde{g} = \tilde{g} (\tilde{I}, \tilde{u},\tilde{v})$ and $ \tilde{f}$ satisfy 
%\footnote{For easier notations we consider the $\sup$ over a domain $D$ also for functions
%that do not depend on the whole set of variables in $D$. 
%From now on we will be consistent with this choice.}
\begin{equation}
%\sup_{ D^d \left(I_0 , \frac{s}{2} \right)  \times D^{2p} \left(0, \frac{r_1}{2} \right)  } 
\sup_{ \mathcal{M}_1 } | \tilde{g} |  \leq  C  \frac{ {(\e M_0)}^2 }{ s {\alpha_1} }
\,\,  ,  \,\,\,\,\,\,\,
\sup_{ \mathcal{M}_1 }  | \tilde{f} |  \leq  { (\e M_0) }^{\nu_1}
\label{stimeteomedia}
\end{equation}
for a suitable $C = C(\sigma, \nu_1)$.
\end{lem}

\nl
Thus, if we set $ \tilde{g} =: \e^2 \bar{g}$ and $\tilde{f} =: \e^{\nu_1}\bar{f} $, 
using (\ref{formaf}) and (\ref{stimeteomedia}) we have
\begin{equation}
\label{H^1_e}
\left\{
\begin{array}{l}
H^1_\e (\tilde{\varphi}, \tilde{I}, \tilde{u}, \tilde{v}) =
                h (\tilde{I}) +  \e \left[ f_0 ( \tilde{I}, \tilde{u}, \tilde{v} ) +
                \e \bar{g}( \tilde{I}, \tilde{u}, \tilde{v} ) \right] +
                \e^{\nu_1} \bar{f} (\tilde{\varphi}, \tilde{I}, \tilde{u}, \tilde{v}) 
\\
\\
f_0 (\tilde{I}, \tilde{u}, \tilde{v}) = f_{00}( \tilde{I} )
                                + \sum_{j=1}^p \Omega_j (\tilde{I})  \frac{\tilde{u}_j^2 + \tilde{v}_j^2}{2}
                                +  O ( {|\tilde{u}, \tilde{v}|}^3 ; \tilde{I})    \,\, .
\end{array} 
\right.
\end{equation}

\nl
From equation (\ref{H^1_e}) we see that the application of averaging theory may cause, in general,  a shift
of order $\e$ of  the elliptic equilibrium, which, before, was in the origin of $\R^{2p}$. 
Therefore, we focus our attention on the Hamiltonian function $f_0 + \e \bar{g}$ with the aim to find a 
real--analytic symplectic transformation restoring the equilibrium in the origin.
An application of the standard Implicit Function Theorem yields the following\footnote{
Since $\Omega_j \neq 0$ for every $j$ in view of (\ref{Ouniformenonrisonanza}), 
we can apply the Implicit function Theorem to obtain, for small enough $\e$, the existence of two functions 
$u_0 = u_0 (\tilde{I} ,\e)$ and  $v_0 = v_0 (\tilde{I}, \e)$ which are
real--analytic for $\tilde{I} \in D^d (I_0 , s/4)$
and such that $\nabla_{\tilde{u}, \tilde{v}} (f_0 + \e \bar{g}) (\tilde{I}, u_0, v_0) = 0$.
Furthermore, using (\ref{stimeteomedia}) together with $\tilde{g} = \e^2 \bar{g}$ and (\ref{sa_1tau_0}), 
one has $u_0, v_0 = O \big( \e  {\left( \log \e^{-1} \right)}^{ 2\tau_0 + 1 } \big)$.
The symplectic transformation in Lemma~\ref{lemmaPhi^2_e} is then generated by
$ x \cdot \tilde{\varphi} + ( p + u_0 (x, \e) ) \cdot ( \tilde{v} - v_0 (x, \e) )$
.}:

\begin{lem}
\label{lemmaPhi^2_e}
Let $\mathcal{M}_2 :=  \T^d_\frac{\s}{7} \times D^d (I_0 , s/4)  \times D^{2p} (0, r_1/4)$;
then, provided $\e$ is sufficiently small,
there exists a (close to the identity)
real--analytic symplectic transformation \begin{equation*}
\Phi^2_\e  :  (x,y,p,q) \in \mathcal{M}_2 
\longrightarrow ( \tilde{\varphi}, \tilde{I}, \tilde{u},\tilde{v} ) \in \mathcal{M}_1 
%:= \T^d_\frac{\s}{6}  \times D^d (I_0 , s/2)  \times D^{2p} (0, r_1/2) 
\end{equation*}
such that $H^2_\e := H^1_\e \circ \Phi^2_\e$ is of the form
\begin{equation*}
H^2_\e (x,y,p,q) =  h(y) + \e \hat{g} (y,p,q) + \e^{\nu_1} \hat{f} (x,y,p,q)
\end{equation*}
with $ \partial_p \hat{g} (y,0,0) = 0 = \partial_q \hat{g} (y,0,0)$,
$\hat{g}$ and $\hat{f}$ real--analytic on $\mathcal{M}_2$.
\end{lem}

\nl
Now, we need to control  the  frequencies associated to the modified Hamiltonian $\hat g(y,0,0)$:

\begin{lem}\label{lem:6}
If $\e$ is small enough then the eigenvalues of the Hamiltonian $\hat{g}$,
i.e. the eigenvalues of\footnote{
$J_{2p}$ denotes the standard $2p \times 2p$ symplectic matrix.}  $J_{2p} \partial^2_{(p,q)} \hat{g} (y,0,0)$,
are given by $2p$ purely imaginary functions $ \pm i \hat{\Omega}_1, \dots , \pm i \hat{\Omega}_p $ verifying
\begin{equation}
\label{hatO-O}
\sup_{y \in D^d (I_0, s/4) }  | \hat{\Omega} (y) - \Omega (y) | = 
									O \left( \e {\left( \log \e^{-1} \right) }^{ 2 \tau_0 + 1}  \right)
\end{equation}
for $\tau_0$ as in (\ref{sa_1tau_0}).
\end{lem}

\nl 
{\bf Proof of Lemma~\ref{lem:6}.}
Consider the quadratic part of $\hat{g}$, 
that is the real--analytic $2p \times 2p$ symmetric matrix $\hat{A} (y) := \partial^2_{(p,q)} \hat{g} (y,0,0)$.
Using the construction of $\Phi^2_\e$ in Lemma~\ref{lemmaPhi^2_e}, 
$\hat{g} = \left( f_0 + \e \bar{g} \right) \circ \Phi^2_\e$, 
estimate (\ref{stimeteomedia}) together with $\tilde{g} = \e^2 \bar{g}$ and the definition of 
$s$ and $\alpha_1$ in (\ref{sa_1tau_0}),
equation (\ref{H^1_e}) for $f_0$  and Cauchy's estimate for derivatives of analytic functions, 
one has
\begin{equation*}
\hat{A} (y) = \diag ( \Omega_1 (y), \dots , \Omega_p (y),  \Omega_1 (y), \dots , \Omega_p (y) ) + 
O  \left( \e { \left(\log \e^{-1} \right) }^{ 2 \tau_0 + 1} \right) \, .
\end{equation*} 
Since $\Omega_j \neq \Omega_k$ for\footnote{From (\ref{Ouniformenonrisonanza}).} $j \neq k$,
an application of the Implicit Function Theorem tells us that the eigenvalues of 
$\hat{g}$ (that a priori might have non--zero real part) 
are $O \left( \e { \left( \log \e^{-1} \right) }^{ 2 \tau_0 + 1} \right)$ close to $\pm i \Omega_j$.
Now, as it is well known, eigenvalues of Hamiltonians always appear in quadruplets $\pm \l, \pm \bar{\l}$;
thus, from the simplicity of the eigenvalues of $\hat{g}$ (holding for $\e$ small enough) one has
that its eigenvalues are purely imaginary as claimed. \hfill$\blacksquare$

\nl
By normal form theory (see corollary 8.7 of \cite{arnold2})
we can find a real--analytic symplectic transformation $O(\e)$--close to the identity
\begin{equation*}
\label{Phi^3_e}
\Phi^3_\e : (\tilde{x} , \tilde{y} , \tilde{p} , \tilde{q}) \in \mathcal{M}_3 := \T^d_\frac{\s}{8} \times 
D^d (I_0, \frac{s}5) \times D^{2p} (0, \frac{r_1}{5})  \longrightarrow  (x,y,p,q) \in \mathcal{M}_2
\end{equation*}
with $\tilde{y} = y$ and such that the transformed Hamiltonian function $H^3_\e := H^2_\e \circ \Phi^3_\e$, 
which is real--analytic on $\mathcal{M}_3$, has the form
\begin{eqnarray*}
\nonumber
H^3_\e ( \tilde{x} , \tilde{y} , \tilde{p} , \tilde{q} ) & = &
h (\tilde{y} ) + \e \hat{g}_0 ( \tilde{y} ) +
					\frac{\e}{2} \sum_{j=1}^{p} \hat{\Omega}_j (\tilde{y}) \left( \tilde{p}_j^2  + \tilde{q}_j^2 \right) 
\\  \nonumber
\\																	
& & \phantom{h (\tilde{y} )} + 		\e \tilde{g}_3 (\tilde{y} ,\tilde{p}, \tilde{q} )  +  
					\e^{\nu_1} \tilde{f}_3 (\tilde{x} , \tilde{y} , \tilde{p} , \tilde{q}) 
%\label{H^3_e}
\end{eqnarray*}	
where $\hat{g}_0 := \hat{g} (\tilde{y},0,0)$,  $\tilde{f}_3 := \bar{f} \circ \Phi^3_\e$ 
and $\tilde{g}_3 := \hat{g}_3 \circ \Phi^3_\e$ verifies
\begin{equation*}
\sup_{ \tilde{y} \in D^d ( I_0, s/5) } | \tilde{g}_3  ( \tilde{y} , \tilde{p} , \tilde{q} ) |
\leq  C {| (\tilde{p} , \tilde{q}) |}^3 \,\,\,  \for  (\tilde{p}, \tilde{q}) \in D^{2p} (0, r_1/5)  \, .
\end{equation*}

\nl
Now let $\tilde{g}_2 ( \tilde{y} , \tilde{p} , \tilde{q} ) 
:= \frac{1}{2} \sum_{i=1}^p \hat{\Omega}_i ( \tilde{p}_i^2 + \tilde{q}_i^2 )$,
we want to put  $ \tilde{g}_2 + \e \tilde{g}_3$ into Birkhoff's normal form up to order $\nu_2$.
In view of inequalities (\ref{Ouniformenonrisonanza}) and (\ref{hatO-O}), provided $\e$ is small enough, we have
\begin{equation}
\label{hatOuninonris}
|  \hat{\Omega} ( \tilde{y} ) \cdot k  | \geq \frac{\alpha_2}{2}  \, , 
\for k \in \Z^p ,\, 0 < {|k|}_1  \leq  \nu_2 , \, \for \tilde{y} \in D^d (I_0, s/5)  \, .
\end{equation}
By Birkhoff's normal form theory\footnote{See, e g., \cite[Theorem 11, p. 43]{HZ} or
 \cite[section 3.4]{miatesi} for a quantitative version.}, one obtains easily the following 

\begin{lem}
\label{lemmaBirkhoff}
If inequality (\ref{hatOuninonris}) is satisfied, 
then there exist $0 < r_\star < r_1^\p \leq r_1/5$ and a real--analytic symplectic diffeomorphism $\Phi^4_\e :  (\theta, r , u, v) \to  ( \tilde{x} , \tilde{y} , \tilde{p} , \tilde{q} )$ mapping 
$$
\label{Phi^4_e}
\mathcal{M}_4 := \T^d_\frac{\s}{8}  \times D^d (I_0 , \frac{s}5)  
\times  D^{2p} (0, r_\star)    \stackrel{\Phi^4_\e}{\to}
\mathcal{M}_3^\p :=   
							\T^d_\frac{\s}{8}  \times D^d (I_0 , \frac{s}5)  \times D^{2p} (0, r_1^\p) 
$$
leaving the origin and the quadratic part of $H^3_\e$ invariant, 
such that  $(\theta, r) = (\tilde{x} , \tilde{y})$ and $H^4_\e := H^3_\e \circ \Phi^4_\e$ is of the form
\begin{eqnarray*}
\label{H^4_e}
\nonumber
H^4_\e ( \theta , r, u, v ) & = & h(r) + \e\hat{g}_0 (r) + 
\frac{\e}{2} \sum_{j=1}^p \hat{\Omega}_j (r) (u_j^2 + v_j^2) + 
\\  \nonumber
\\
& & +	\	\e\, Q_\star (r,u,v)  +  \e R_\star (r,u,v) +  \e^{\nu_1} \tilde{f}_4 ( \theta , r, u, v )  	
\end{eqnarray*}
where:
\begin{itemize}
	\item $Q_\star$ is a polynomial of degree $\left[ \frac{\nu_2}{2} \right]$ in the variables $I = (I_1, \dots, I_p) $ 
	having the form
	\begin{equation*}
	\langle  \hat{\Omega} (r), I \rangle + \frac{1}{2} \langle T (r) I, I  \rangle + \cdots 
	\qquad \mbox{with}  \qquad I_j := \frac{1}{2} (u_j^2 + v_j^2) \, 
	\end{equation*}
	with $T(r)$ a $2p \times 2p$ real--analytic matrix;
	
	\item $R_\star$ is a real--analytic function verifying $ | R_\star (r,u,v) | \leq C {|(u,v)|}^{\nu_2+1}$
	for every $(u,v) \in D^{2p} (0, r_\star)$ and $r \in D^d (I_0, s/5)$;
	
	\item $\tilde{f}_4 := \tilde{f}_3 \circ \Phi^4_\e$ is real--analytic on $\mathcal{M}_4$.
\end{itemize}
\end{lem}

\nl
We may conclude the proof of Proposition~\ref{teoremaconiugazione}.
Following \cite[pp. 1561--1562]{fejoz}, 
we pass to symplectic polar coordinates in order to move $R_\star$ to the perturbation
of $H^4_\e$ with the help of a rescaling by a factor $\e$. 
Let $\rho^0 = (\rho^0_1, \dots, \rho^0_p)$ in $ {(\R_+) }^p$ be sufficiently close to the origin;
consider, for a suitable ${\sigma_\star} > 0$, the real--analytic symplectic transformation
$\Phi^5_\e : (\theta, r, \zeta, \rho) \to(\theta, I_0 + r, z)$ mapping 
\begin{eqnarray*}
\mathcal{M}_5 := \T^d_\frac{\s}{8} \times D^d (0, \frac{s}5) \times
\T^p_{\sigma_\star} \times D^p \big( 0, \frac{|\rho^0|}2 \big)  
\stackrel{\Phi^5_\e}{
\longrightarrow}
\mathcal{M}_4
\end{eqnarray*}
where 
\begin{equation}
\label{coordinatepolari}
z_j = u_j + i v_j := \sqrt{2 \left( \rho^0_j + \rho_j \right) } \, e^{- i \zeta_j} \, .
\end{equation}
The transformed Hamiltonian function $H^5_\e := H^4_\e \circ \Phi^5_\e$, 
real--analytic on $\mathcal{M}_5$, assumes the form
\begin{eqnarray*}
\nonumber
&&\!\!\!H^5_\e  ( \theta, r, \zeta, \rho )  =  h(I_0 + r) +  \e \hat{g}_0 (I_0 + r) + 
			\frac{\e}{2} \sum_{j=1}^p 	\hat{\Omega}_j (I_0 + r) (\rho^0_j + \rho^0_j)
			\\
			\\
& &\phantom{aAAAAAA} 	+\  \e	Q_{ I_0 + r} (\rho^0 + \rho) +
			\e R (I_0 + r , \zeta, \rho^0 + \rho) + \e^{\nu_1} \tilde{f}_5 ( \theta, r, \zeta, \rho; \rho^0 )
\end{eqnarray*}
where 
\begin{itemize}

	\item $Q_{I_0 + r} := Q_\star \circ \Phi^5_\e$ is a polynomial of degree $\left[ \frac{\nu_2}{2} \right]$ with respect to 
	$\rho^0 + \rho$, depending also on $I_0 + r$;
	
	\item $R := R_\star \circ \Phi^5_\e$ verifies  
	\begin{equation*} 
	|R (I_0 + r , \zeta, \rho^0 + \rho)| \leq C {\left| \rho^0 \right|}^\frac{\nu_2+1}{2}
	\end{equation*}
	for every $ \rho \in D^{2p} \left( 0, |\rho^0| /2 \right)$, $r \in D^d (0, s/5)$ and $\zeta \in \T^p_{\sigma_\star}$;
	
	\item $\tilde{f}_5 := \tilde{f}_4 \circ \Phi^5_\e$ is real--analytic on $\mathcal{M}_5$.
	
\end{itemize}	 
Now, let $A_\e$ be the homothety given by
\begin{equation*}
A_\e :  ( \theta, r, \zeta, \rho )   \longrightarrow (\theta, \e r , \zeta, \e \rho) \, .
\end{equation*}
Even though $A_\e$ is not a symplectic map it preserves the structure of Hamilton's equations if we consider
the Hamiltonian function $H^6_\e := \frac{1}{\e} H^5_\e \circ A_\e$. 
Explicitly we have
\begin{eqnarray}
 H^6_\e ( \theta, r, \zeta, \rho  )  &=&  \frac{1}{\e}  h (I_0 + \e r) + \hat{g}_0 (I_0 +\e r) +
				 \frac{1}{2} \hat{\Omega} (I_0 + \e r) \cdot (\rho^0 + \e \rho) 
\nonumber				
\\  \nonumber
\\
& & +   Q_{\e, I_0 + \e r} ( \rho^0 + \e \rho) +	   	
	R (I_0 + \e r, \e \rho , \zeta ; \rho^0 )	\nonumber \\  \nonumber\\ 
	&&  + \e^{\nu_1 - 1} \tilde{f}_6 ( \theta, r, \zeta, \rho; \rho^0)	\, .
\label{H^6_e}
\end{eqnarray}
where $\tilde{f}_6 := \tilde{f}_5 \circ A_\e$.
Now we fix $\rho^0 \in {(\R_+)}^p$ with $| \rho^0 | = 2 \e$ so that  $| R | \leq C \e^\frac{\nu_2 + 1}{2}$.
Thus, if we choose $\nu_1$ and $\nu_2$ so that
\begin{equation}
\label{eqN_1N_2}
\nu_1 - 1 = \left[ \frac{\nu_2 + 1}{2} \right] := \nu
\end{equation} 
we may write
\begin{equation}
\label{perturbazionehatH^6_e}
R (I_0 + \e r, \e \rho , \zeta ; \rho^0 )	+ \e^{\nu_1 - 1} \tilde{f}_6 ( \theta, r, \zeta, \rho ) =:	
\e^\nu  P_\e (\theta, r, \zeta, \rho) 
\end{equation}
for a suitable function $P_\e$ real--analytic on $ \T^d_\frac{\s}{8} \times D^d (0, s/5) \times
\T^p_{\sigma_\star} \times D^p ( 0, \e ) $ 
.

\nl
We have proved Proposition~\ref{teoremaconiugazione} with $\hat{H}_\e = H^6_\e$
as in (\ref{H^6_e}), (\ref{perturbazionehatH^6_e}).  \hfill$\blacksquare$ 

\medskip

\subsubsection{Step 2: Amounts of nondegeneracy of the normal form}

\begin{pro} 						
\label{teoremanondeg}
Let $N_\e$ be as in (\ref{Neps}).
If $\e$ is small enough, the frequency map
\begin{equation*}
\hat{\Psi}_\e : (r,\rho) \in B^d (0, s/5) \times B^p (0, \e) \longrightarrow 
														\left( \frac{\partial}{\partial r} N_\e , %(r_0 + \e r,\rho^0 + \e \rho)  
														\frac{\partial}{\partial \rho} N_\e  %(r_0 + \e r,\rho^0+ \e \rho)   
														\right)
\end{equation*}
is R--nondegenerate.

\nl
Moreover, let $\bar{\mu}$ and $\bar{\beta}$ denote respectively the index and the amount of nondegeneracy of the 
unperturbed frequency map (\ref{applicazionefreq}) with respect to a closed ball 
$\bar{B}^d (I_0, t) \subset \B$, for some $t > 0$ independent of $\e$.
Then, if we define $\K_\e := \bar{B}^d (0, s/10) \times \bar{B}^p (0,\e/2)$ and let 
$\hat{\mu}_\e$ denote the index of nondegeneracy of $\hat{\Psi}_\e$ 
with respect to $\K_\e$ and 
\begin{equation*}
\hat{\beta}_\e :=  \min_{ c \in \mS^{d+p-1} } \min_{ (r,\rho) \in \K_\e } 
					\max_{0 \leq \mu \leq \bar{\mu}}  \left|  \partial^\mu_{(r, \rho)}
										{ | \langle c, \hat{\Psi}_\e \rangle | }^2  \right|  \, ,
\end{equation*} 
one has
\begin{equation}
\label{hatmu_e}
\hat{\mu}_\e \leq \bar{\mu}  \qquad \mbox{and} 
						\qquad \hat{\beta}_\e \geq \frac{\e^{\bar{\mu}+2} \bar{\beta} }{8} \, .
\end{equation}
\end{pro}

\nl
{\bf Proof of Proposition \ref{teoremanondeg}.}
From (\ref{Neps}), it follows that the frequency map of $N_\e$ is given by 
\begin{eqnarray*}
& & \hat{\Psi}_\e (r,\rho)  = \left(  \o (I_0 + \e r) + O(\e) ,  \frac{\e}{2} \hat{\Omega} (I_0 + \e r)	+ O(\e^2)\right)	
\end{eqnarray*}
and it is real--analytic on  $D^d (0, s/5) \times D^p ( 0, \e )$.
Using (\ref{hatO-O}) one has
\begin{equation}
\label{formulahatPsi_e}
\hat{\Psi}_\e (r,\rho) = \left(  \o (I_0 + \e r) + O(\e), \, \frac{\e}{2} \big( \Omega (I_0 + \e r)	+ O(\e) \big) \right)	\,.
\end{equation}

\noindent
Now,  let $\bar{\mu} \in \N_+$ and $\bar{\beta} > 0$ 
denote respectively the index and the amount of nondegeneracy
of $\Psi := ( \o, \Omega)$ with respect to $\bar{B}^d (I_0 , t)$, for some positive $t$ independent of $\e$.
Set
\begin{equation}
\label{formulaPsi_0}
\Psi_0 (r) := \left(  \o (I_0 + \e r) ,   \Omega (I_0 + \e r) \right)	\, ,
\end{equation}
$\K_0 := \bar{B}^d (0, s/10)$
and use definition \ref{defindexamount} to get
\begin{equation*}
\min_{ r \in  \K_0 } \max_{0 \leq \mu \leq \bar{\mu} }  
							\left| \partial^\mu_r {\left| \langle c, \Psi_0 (r)  \rangle \right|}^2 \right| 
							\geq \e^{ \bar{\mu} } \bar{\beta} > 0
\end{equation*}
for every $c \in \mS^{d+p-1}$.

\nl
Next, denote by $\Psi_\e$ the real--analytic function over $ D^d (0, s/5) \times D^p ( 0, \e )$
obtained multiplying the last $p$ component of $\hat{\Psi}_\e$ by a factor $2/\e$. 
Then, observe that equations (\ref{formulahatPsi_e}) and (\ref{formulaPsi_0}) 
imply $\Psi_\e (r, \rho) = \Psi_0 (r) + O(\e)$.
Therefore, denoting $\K_1 := \bar{B}^p (0 , \e/2)$ and assuming $\e$ small enough, 
one has
\begin{equation*}
\beta_\e := \min_{ (r,\rho) \in \K_0 \times \K_1 } 
							\max_{0 \leq \mu \leq \bar{\mu} }  
							\left| \partial^\mu_{(r ,\rho)}  {\left| \langle c, \Psi_\e (r, \rho)  \rangle \right|}^2 \right| 
							\geq \frac{\e^{\bar{\mu}} \bar{\beta}}{2} > 0
\end{equation*}
for every $c \in \mS^{d+p-1}$.

\nl
Now, if we write $\Psi_\e = ( \Psi_\e^{(1)} , \Psi^{(2)} ) \in \R^d \times \R^p$,
from what observed before, it results
\begin{equation*}
\hat{\Psi}_\e = \left( \Psi_\e^{(1)} , \frac{\e}{2} \Psi_\e^{(2)} \right) \, .
\end{equation*}
Define for $c = (c_1, c_2) \in \R^d \times \R^p $ with ${|c|} = 1$ the function
\begin{eqnarray*}
f ( r, \rho, c_1, c_2 )& :=&
							\max_{0 \leq \mu \leq \bar{\mu} }  
							\left| \partial^\mu_{(r, \rho)} {\left| \langle c, \Psi_\e  \rangle \right|}^2 \right|  
						\\
						&	=&
							\max_{0 \leq \mu \leq \bar{\mu} }  
							\left| \partial^\mu_{(r, \rho)} {\left| \langle c_1, \Psi^{(1)}_\e  \rangle 
																		+  \langle c_2, \Psi^{(2)}_\e  \rangle  \right|}^2 \right|  \,\, ;
\end{eqnarray*}				
furthermore set 
\begin{equation*}
t_\e :=  \sqrt{ {|c_1|}^2 + \frac{\e^2}{4} {|c_2|}^2 }
\end{equation*}
and $\bar{c}_1 = c_1 t_\e^{-1}$, $\bar{c}_2 = \e c_2 {(2 t_\e)}^{-1}$ so that ${|(\bar{c}_1, \bar{c}_2)|} = 1$.
Then one has
\begin{eqnarray*}
\max_{0 \leq \mu \leq \bar{\mu} }  
							\left| \partial^\mu_{(r, \rho)} {\left| \langle c, \hat{\Psi}_\e  \rangle \right|}^2 \right|
							& = & f \left(r, \rho, c_1, \frac{\e}{2} c_2 \right)
							= t_\e^2  f \left(r, \rho, \frac{c_1}{t_\e}, \frac{\e}{2} \frac{c_2}{t_\e} \right) 
							\\
							\\
							& \geq & \frac{\e^2}{4} f \left(r, \rho, \bar{c}_1, \bar{c}_1 \right) 
							\geq \frac{ \e^{\bar{\mu} + 2} \bar{\beta} }{8} > 0  \, 
\end{eqnarray*}
and it follows immediately
\begin{equation*}
\min_{ (r,\rho) \in  \K_0 \times \K_1} \max_{0 \leq \mu \leq \bar{\mu}}  \left|  \partial^\mu_{(r, \rho)}
										{ | \langle c, \hat{\Psi}_\e \rangle |}^2   \right| 
										\geq \frac{\e^{\bar{\mu}+2} \bar{\beta} }{8} > 0
\end{equation*}
for every $c \in \mS^{d+p-1}$.
Since $\K_0 \times \K_1 = \bar{B}^d (0, s/10) \times \bar{B}^p (0, \e/2) = \K_\e$ 
we have verified (\ref{hatmu_e}).
In view of remark \ref{ossRnondeg} we also conclude that $\hat{\Psi}_\e$ is R--nondegenerate 
on $B^d (0, s/5) \times B^p (0, \e)$, provided that $\e$ is small enough.

\nl
Proposition~\ref{teoremanondeg} is proved.  \hfill$\blacksquare$

\subsubsection{Conclusion of the Proof of Theorem \ref{teoremadegenere}}
\label{subsec:conc}
We want to apply \Russ's Theorem \ref{teorussmassimale} to the properly degenerate case of $\hat{H}_\e$ in (\ref{Hconiugata}).
With Propositions \ref{teoremaconiugazione} and \ref{teoremanondeg} we are in a position to meet 
the hypothesis of R--nondegeneracy of the frequency application required in Theorem \ref{teorussmassimale}.
However, the ``degenerate'' case of $ \hat{H}_\e$ requires
that the size of its perturbation is of a sufficiently small order in $\e$.
From (\ref{Hconiugata}) we see that the size of the perturbation of $\hat{H}_\e$ is order $\e^\nu$
where $\nu$ can be chosen to be arbitrarily big\footnote{ 
Recall (\ref{eqN_1N_2}) and the fact that both $\nu_1$ and $\nu_2$ can be arbitrarily fixed at the beginning of
the process described in section \ref{sezdimlemma}.
}
but independent of $\e$.

\nl
Next we provide an explicit expression for the admissible size of the perturbation in \Russ's Theorem,
i.e. $\e_0$ in (\ref{tagliaperturbazione}).

\begin{lem}[\bf \Russ]
\label{lemmae_0}
Let $H, \mathcal{Y}$, $\e^\star$, $\mathcal{A}$ and $\tau$ be as in Theorem \ref{teorussmassimale} 
and let $\o := \nabla h$ be R--nondegenerate (as in the hypotheses of Theorem \ref{teorussmassimale}).
Consider the following quantities:

\begin{enumerate}

				\item \label{muebeta}
								Let $\K \subset \Y$ be any chosen compact set;
        				let $\mu$ be any integer greater than the 
        				index of nondegeneracy of $\o$ with respect to $\K$ 
        				and let $\beta$ be the ``amount of nondegeneracy'' corresponding to $\mu$.

        \item \label{costantiK}
        				Let $\vartheta \in (0,1)$ be chosen such that\footnote{
        				Recall definition (\ref{notazioneinsiemi})} 
        				$ \T^d_\vartheta \times \left( \K + 4\vartheta \right) \subset \A$
        			  and define $C_1 := {|\o|}_{\K + 3\vartheta}$.
                Let $d_0$ be the diameter of $\K$, i.e. $d_0 := \sup_{x,y \in \K} |x-y|$.

        \item \label{PhieT_0}
        				%Let $\Phi : [1, \infty) \longrightarrow \R$ be defined as
 								%\begin{equation}
 								%\label{sceltaPhi}
 								%\Phi (T) := T^{-\tau} \qquad (\mbox{with} \,\,\,  \tau > n \bar{\mu})  \, .
 								%\end{equation}
 								%Furthermore 
 								Let $ T_0 \geq  e^\frac{n+1}{\tau}$ such that the following inequality holds
                \begin{equation}
                \label{condizioneT_0teoruss}
                \int_{T_0}^\infty   \frac{\log T}{T^2} \, dT  \leq  \vartheta \, .
                \end{equation}

        \item \label{get_0}
        				Define 
                \begin{eqnarray}
                %\label{M^starmassimale}
                %M^\star & := & 2^{ \mu + 2} \mu ! \vartheta^{- \mu } (C_1 + 1)  \,  ,
                %\\ \nonumber
                %\\
                \label{C^starmassimale}
                C^\star & := & 2^{\mu + 1} 
                	\frac{ {(\mu + 1)}^{ \mu + 2} }{\vartheta^{\mu + 1}} (C_1 + 1) 
                \end{eqnarray}
             		and set 
                \begin{eqnarray}
                \label{defgmassimale}
                \g  & := & \left( d_0^n  C^\star  \right)^{- \frac{ \mu }{2} } 
                															\beta^{\frac{ \mu + 1 }{2} }  
                															{ \e^\star }^{\frac{ \mu }{2} }  \, ,
								\\ \nonumber
								\\                															
                \label{deft_0massimale}
                t_0 & := & \frac{ \g T_0^{-(\tau + n + 1)} \vartheta }{C_1 + 1}   \,  .
                \end{eqnarray}
                
				\item Finally set 
								\begin{eqnarray}
        				\nonumber
								E_1 & := & \g  T_0^{-(\tau + n + 1)} \vartheta
								\\      \nonumber
								\\      
								E_2 & := & \frac{ \beta  t_0^{ \mu} }{ T_0 C^\star }\, .
								\label{defE_imassimale}
								\end{eqnarray}        
\end{enumerate}

Then $\e_0$ in (\ref{tagliaperturbazione}) can be taken to be
\begin{equation}																																						
\label{e_0massimale}
\e_0 : = c_0 \frac{\vartheta}{C_1} { \left( \min \{ E_1 , E_2 \}  \right) }^2
\end{equation}
for a suitable $c_0 = c_0 ( n, \mu )$.
\end{lem}

\begin{oss}\rm 
\label{ossPhi}
The above result follows from   \cite{russmann}  by considering the case of maximal tori only\footnote{Compare in particular the estimates listed
on page 171 of \cite{russmann};  see, also,  chapter 2 of \cite{miatesi}.}. More precisely:

\begin{enumerate}

\item[(i)] 
		The maximal case corresponds to the easier case $p=q=0$ in \cite{russmann}.
		Notice, however, that it is not sufficient to substitute the values $p=q=0$ in \Russ's
		estimates (as ,when $p=0$ for instance, many terms in \cite[p. 171]{russmann} become meaningless)
		but, rather, one has to go through the most of Theorem 18.5 in \cite{russmann} to get the value
		of $C^\star$ and $\g$ in (\ref{C^starmassimale}) and (\ref{defgmassimale}) 
		and through the first part of Lemma~13.4 in \cite[pp. 158--161]{russmann} to
		get the value of $t_0$ in (\ref{deft_0massimale}).

		\nl 
\item[(ii)]
{\rm (Control of small divisors)}
In \cite[section 1.4]{russmann} \Russ \, introduces a so called ``approximation function'' $\Phi$
in order to control the small divisors.
Our choice is to take $\Phi (T) = T^{-\tau}$ with $\tau > n \mu$.
Comparing \cite[section 1.4]{russmann}, one sees that such $\Phi$ does not verify property 3, 
i.e. $T^\l \Phi (T) \stackrel{T \rightarrow \infty}{\longrightarrow} 0$ for any $\l \geq 0$.
However, when we consider $H = \hat{H}_\e$, we will see below that one has $T_0 = T_{0,\e} = O (\e^{-2})$;
then, equation 14.10.10 together with 13.1.4 and inequality 14.10.11 in \cite{russmann}
would cause $O(\e^\nu)$ to be an inadmissible size for a perturbation.
Nevertheless we claim that the only decay property which is actually needed 
in \Russ's Theorem \ref{teorussmassimale} is
\begin{equation*}
\lim_{T \rightarrow \infty} T^\l \Phi (T) = 0 \qquad \mbox{for \, all} \qquad 0 \leq \l < n \mu
\end{equation*}
so that our choice is perfectly suitable.
\end{enumerate}

\end{oss}\rm

\nl
Now we are going to analyze what happens to the estimate in Lemma~\ref{lemmae_0} 
when we consider $\hat{H}_\e$ as Hamiltonian function.
In particular we are going to show that each one of the quantities appearing in 
Lemma~\ref{lemmae_0} 
can be controlled by constants involving initial parameters 
related only to $H_\e$ in (\ref{Hamiltonianainiziale}) and (\ref{formaf})
times powers of $\e$.

\nl
We point out that in the application of Theorem \ref{teorussmassimale}
with $H = \hat{H}_\e$ in (\ref{Hconiugata}) 
we have the following correspondences\footnote{
See Theorem \ref{teorussmassimale}, Proposition~\ref{teoremaconiugazione},
(\ref{notazioneinsiemi}) and (\ref{sa_1tau_0}) for notations.}:

\begin{equation}
\label{corrispondenze}
\begin{array}{l}
\T^n = \T^d \times \T^p  \,\,\,  ,  \,\,\,\,\,  x = (\theta,\zeta)  
\\
\\
\Y = \Y_\e :=  B^d (0, s/5) \times  B^p (0, \e)  \,\,\,  ,  \,\,\,\,\,  y = (r,\rho)
\\
\\
\A = \A_\e :=  \T^d_{\frac{\s}{8}}  \times \T^p_{\sigma_\star}  \times D^d (0, s/5 )  \times D^p (0 ,\e)
\\
\\
P = \e^\nu P_\e  \,\,\,  ,  \,\,\,\,\, N = N_\e    \,\, .
\end{array}
\end{equation}

\nl
Accordingly to \ref{lemmae_0}.\ref{muebeta} 
we consider the frequency application of the integrable part of $\hat{H}_\e$, that is
$\hat{\Psi}_\e (r,\rho)$ as in Proposition~\ref{teoremanondeg}.
We already proved that $\hat{\Psi}_\e$ is R--nondegenerate for 
$ (r ,\rho) \in B^d (0 ,s/5) \times B^{2p} (0, \e)$. 
Now, in view of \ref{lemmae_0}.\ref{muebeta} and the correspondences in (\ref{corrispondenze})
we need to fix a compact set $\K = \K_\e \subset \A_\e$.
For our convenience we take $\K_\e := \bar{B}^d (0, s/10) \times  \bar{B}^p (0, \e/2) $
so that the first inequality in (\ref{hatmu_e}) 
allows us to consider\footnote{
Recall Proposition \ref{teoremanondeg} for the definition of $\bar{\mu}$.}
$\mu = \bar{\mu}$ 
as an integer greater than the actual index of nondegeneracy of $\hat{\Psi}_\e$ with respect to $\K_\e$.
Also, in view of the second inequality in (\ref{hatmu_e}), 
we can take
\begin{equation}
\label{beta_e}
\beta = \beta_\e := \frac{\e^{\bar{\mu} + 2} \bar{\beta}}{8}
\end{equation}
in (\ref{defgmassimale}) and (\ref{defE_imassimale}).

\nl
Next, we choose $\vartheta = \vartheta_\e := \e/16$
so that, for $\e$ sufficiently small and in view of (\ref{ordines}),
one has $\K_\e + 4 \vartheta_\e \subset \A_\e$ as required in \ref{lemmae_0}.\ref{costantiK}.
Accordingly to Theorem \ref{teorussmassimale} 
we also need to fix a positive number $ \e^\star < \meas_{d+p}\, \K_\e$. 
In view of our definition of $\K_\e$ and (\ref{ordines}) a suitable choice is given by
\begin{equation}
\label{e^star_e}
\e^\star = \e^{p+1}
\end{equation}
for $\e$ small enough.

\nl
Now, observe that the quantities $C_1$ and $d_0$ in Lemma~\ref{lemmae_0}, point \ref{costantiK}
do not cause any change in the order in $\e$ of the size of the admissible perturbation. 
In fact, using equation (\ref{formulahatPsi_e}) and taking $\e$ sufficiently small, we have
\begin{eqnarray*}
C_1 = C_{1, \e} & := & { |\hat{\Psi}_\e| }_{ \K_\e + 3 \vartheta_\e  }  \leq  
\\
\\
& \leq &	\sup_{ r \in D^d (0, s/5 ) } { | \o ( I_0 + \e r) | } + 
						\e \sup_{ r \in D^d (0, s/5 ) }  { | \Omega ( I_0 + \e r) | }  + O (\e) \leq
\\
\\
& \leq & \sup_{ r \in D^d (I_0 , s/5 ) }  { | \o ( r ) | } + 
						\e \sup_{ r \in D^d (I_0 , s/5 ) }  { | \Omega ( r ) | }  + O (\e)  \leq	M_1
\end{eqnarray*}
where $M_1$ is defined in (\ref{normaoO}).
Since the estimate for $\e_0$ is decreasing with respect to $C_1$,
we can substitute $C_1$ in %(\ref{M^starmassimale}), 
(\ref{C^starmassimale}) and (\ref{deft_0massimale}) with $M_1$. 
The estimate for $\e_0$ is also decreasing with respect to $d_0$ so that 
when we consider $\K = \K_\e$ we may simply replace $d_0$ by $1$.

\nl
Let us now analyze the quantities in \ref{lemmae_0}.\ref{PhieT_0} and \ref{lemmae_0}.\ref{get_0}.
%\footnote{The different choice of $\Phi$ with respect to \Russ's ''approximation function'' has been already
%explained briefly after Theorem \ref{teorussmassimale}. A fully detailed explanation can be found in 
%\cite{miatesi}[section 4.3.5].}
First of all observe that in view of (\ref{hatmu_e}) and $ n = d+p$
we can fix a priopri an exponent $\tau \geq (d+p) \bar{\mu}$ 
satisfying the requirement in \ref{lemmae_0}.\ref{PhieT_0}.
Furthermore, given the previous choice of $\vartheta_\e$, 
inequality (\ref{condizioneT_0teoruss}) becomes
\begin{equation*}
\int_{T_0 }^\infty     \frac{\log T}{T^2} \, dT  \leq \frac{\e}{16} 
\end{equation*}
which can be easily fulfilled, together with $T_0 \geq e^{ \frac{d+p+1}{\tau} }$,  by choosing
\begin{equation}
\label{T_0e}
T_0 = T_{0 , \e} :=  \frac{1}{\e^2}
\end{equation}
for $\e$ sufficiently small.
For what concerns the quantities defined in \ref{lemmae_0}.\ref{get_0} we see that
since $\vartheta = \vartheta_\e := \e/16$ and the estimate for $\e_0$ is decreasing in 
%both 
$C^\star$,
%and $M^\star$,
we can choose
\begin{equation}
\label{C^star_e}
%M^\star_\e := c \sup_{B_{r_0}} |(\o, \Omega)| \e^{ - \bar{\mu}}  \qquad   \mbox{and}  \qquad   
C^\star = C^\star_\e :=  2^{5 (\bar{\mu} + 1) } { (\bar{\mu} + 1) }^{ \bar{\mu} + 2} (M_1 + 1) \e^{ - (\bar{\mu}+1)}
\end{equation}
having also used $C_1 = C_{1, \e} \leq M_1$.
From the fact that we can replace $d_0$ by $1$ 
together with equations (\ref{beta_e}), (\ref{e^star_e}) and (\ref{C^star_e}),
one has
\begin{equation}
\label{gamma_e}
\g = \g_\e := c_1 {\left( M_1 + 1 \right)}^{ - \frac{ \bar{\mu} }{2} } 
											\e^{ { (\bar{\mu} + 1)}^2 + \frac{(p+1) \bar{\mu} }{2} } 
											\bar{\beta}^{\frac{ \bar{\mu}+1}{2} }  { \e^\star }^{\frac{\bar{\mu}}{2} }   
\end{equation}
for a suitable constant $c_1 < 1$ depending only on $\bar{\mu}$.
Moreover, given once again the previous choice of $\vartheta = \vartheta_\e$ together with
equation (\ref{T_0e}) and the above definition of $\g_\e$,
we can replace $t_0$ in (\ref{deft_0massimale}) by
\begin{equation}
\label{t_0e}
t_0 = t_{0, \e} = c_1  { \left( M_1 + 1 \right) }^{- \frac{ \bar{\mu} }{2} }
						\e^{ { (\bar{\mu} + 1)}^2 + \frac{(p+1) \bar{\mu}}{2} + 2( \tau + d + p) + 3 }
\end{equation}
for $c_1 < 1$ as above.

\nl
From (\ref{defE_imassimale}) we see that $E_1$ and $E_2$ have simple polynomial dependence on the quantities 
$\g, T_0^{-1}, \vartheta, \beta$ and $t_0$. 
Our previous analysis shows that when we consider $\hat{H}_\e$ as Hamiltonian function, 
these quantities can be replaced respectively by\footnote{
See (\ref{gamma_e}), (\ref{T_0e}), (\ref{beta_e}), (\ref{t_0e}) 
and recall $\vartheta = \vartheta_\e := \e/16$.
}

\begin{eqnarray*}
& & \g_\e = O \left( \e^{ { (\bar{\mu} + 1)}^2 + \frac{(p+1) \bar{\mu}}{2} } \right) \,\, ,
\qquad
T^{-1}_{0, \e} = O (\e^2)  \,\, ,  \qquad \vartheta_\e = O(\e)   \\
\\
& &  \beta_\e = O \left( \e^{ \bar{\mu} + 2 } \right)  \,\,  ,
\qquad t_{0,\e} = O \left( \e^{ { (\bar{\mu} + 1)}^2 + \frac{ (p+1) \bar{\mu} }{2} + 2( \tau + d + p) + 3 } \right) \, .
%{ \left( M^\star_\e \right) }^{-1} = O (\e^{\bar{\mu}})
\end{eqnarray*}
Therefore, in view of (\ref{e_0massimale}) the size of the perturbation allowed by \Russ's Theorem when we consider 
$H = \hat{H}_\e$, is order $\e^{\nu_0}$ with\footnote{
Using (\ref{N_0}) we are able to define the values of $\nu_1$ in (\ref{definizioneK_1}) and
(\ref{ouniformenonrisonanza}) and $\nu_2$ in (\ref{Ouniformenonrisonanza}) through equation (\ref{eqN_1N_2}). 
}
\begin{equation}
\label{N_0}
\nu_0 :=  2 {\bar{\mu}}^3   +   (p+5) {\bar{\mu}}^2   +   [ 14 + 4 (\tau + d + p) ] \bar{\mu} + 13    \, .
\end{equation} 
In particular we have a condition of the form $\e_0 \leq \bar{c} \e^{\nu_0}$ 
where $\bar{c}$ is some positive constant independent of $\e$ 
and depending only on quantities related to the initial Hamiltonian $H_\e$,
namely $\bar{\mu}, \bar{\beta}$ and $\bar{\K}$ as in Proposition~\ref{teoremanondeg},
the Diophantine constant $\tau \geq (d+p) \bar{\mu}$
and $M_1$ as in (\ref{normaoO}) with $r_0$ as in (\ref{dominiooloH}).
By Proposition~\ref{teoremaconiugazione} we know that we can assume 
the size of the perturbation of $\hat{H}_\e$ to 
be order $\e^\nu$ for any fixed integer $\nu \geq 4$ independent of $\e$. 
Thus, by simply taking\footnote{Notice that this can be done since $\nu_0$ only depends 
on $d,p,\tau$ and $\bar{\mu}$.}  $\nu > \nu_0$, 
we can apply \Russ's Theorem to $\hat{H}_\e$ and obtain Theorem \ref{teoremadegenere} as a consequence. \hfill$\blacksquare$

\section{Proof of Theorem~\ref{thm.an}}
\label{sezfinale}
As it follows from  the analysis described in \cite[Sect. 6, pp. 1563--1569]{fejoz},
the motions of $(n+1)$ bodies (point masses)  interacting only through gravitational attraction,  restricted to the invariant symplectic submanifold of vanishing total linear momentum, are governed by the real--analytic Hamiltonian
\begin{equation}\label{ham.n}
F= H^0(\Lambda) + \e  \Big(H^1(\Lambda,\xi,\eta,q,p) + H^2(\l,\Lambda,\xi,\eta,q,p) \Big)
\end{equation}
where:
\begin{itemize}
\item[(i)] $
\left( \l , \Lambda , \xi , \eta , q , p  \right) \in  \T^n \times {(0,\infty) }^n 
\times \R^n \times \R^n \times \R^n \times \R^n$
are standard symplectic coordinates;

\item[(ii)] $\Lambda_j = \mu_j \sqrt{M_j \, a_j}$, where $a_j>0$ are 
the  semi major--axis of the ``instantaneous'' Keplerian ellipse formed by the ``Sun'' (major body)   and the $j^{\rm th}$ ``planet'', while   
$$\frac{1}{\e \mu_j} = \frac{1}{m_0} + \frac{1}{\e m_j}\ ,\qquad
M_j:= m_0 + \e m_j\ ,
$$
$m_0$ and $\e m_j$ being, respectively the mass of the Sun and the mass of the $j^{\rm th}$--planet; 

\item[(iii)] the phase space $\M$  is the  open subset of
$\T^n \times {(0,\infty) }^n 
\times \R^n \times \R^n \times \R^n \times \R^n$ subject to the collisionless constrain
$$0 < a_n < a_{n-1} < \dots < a_1$$
and endowed with the standard symplectic form $\displaystyle\sum_{j=1}^n d\l_j \wedge  d\Lambda_j+ 
d\xi_j \wedge  d\eta_j + dq_j \wedge  dp_j$;

\item[(iv)] $H^0:=F_{\rm Kep}$ is the Keplerian integrable limit given by 
$$
H^0:=F_\Kep := \sum_{j=1}^n -\frac{\mu_j^3 M_j^2}{2 \Lambda_j^2}\ ,
$$
describing $n$ decoupled  two--body systems formed by the Sun and the $j^{\rm th}$ planet;

\item[(v)] the ``secular'' Hamiltonian $H^1$ has the form\footnote{
There is a difference of a factor $\frac{1}{2}$ with the notations used in Ref. \cite{fejoz}. The computations are performed in \cite{laskar}.}

\begin{equation}\label{H1}
H^1=C_0 +\sum_{j=1}^n \s_j \frac{\xi_j^2+\eta_j^2}{2} + \sum_{j=1}^n \varsigma_j \frac{q_j^2+p_j^2}{2} + O(4)
\end{equation}
where $C_0$, $\s_j$ and $\varsigma_j$ depend on $\Lambda$; ``$O(4)$'' denotes terms of order greater than or equal to four in $(\xi,\eta,q,p)$ (and depending on $\Lambda$);

\item[(vi)]
$H^2$ has vanishing average over  $\l\in\T^n$; $H^i$ depend also (in a regular and non influential way) on $\e$.

\end{itemize}

\begin{oss}\rm \rm The variables $(\l , \Lambda , \xi , \eta , q , p)$ are obtained from standard Poincar\'e variables after a rotation in  $(\xi , \eta , q , p)$ needed to diagonalize the quadratic part of the secular Hamiltonian; the ``eigenvalues'' $\sigma_j$ and $\varsigma_j$ are the the first Birkhoff invariants of the secular Hamiltonian; compare \cite[pp. 1568, 1569]{fejoz}. \end{oss}\rm 

\nl
The frequency map of the planetary Hamiltonian $F$ is given by 
$$ \{\nu_1 , \dots , \nu_n  , \s_1 , \dots , \s_n, \vs_1 , \dots , \vs_n \}$$
where the $\nu_j$'s are the Keplerian frequencies 
\begin{equation} \label{movimentimedi}
\nu_j := \frac{\partial F_\Kep}{ \partial \Lambda_j} = 
												\frac{\sqrt{M_j}}{a_j^\frac{3}{2}} = \frac{\mu_j^3 M_j^2}{\Lambda^3_j} \, .
\end{equation}
It is customary  to consider the frequency map as a function of the semi--major
axes $a$ (rather than of the actions $\Lambda$); we shall therefore call the  
``planetary frequency map'' the application\footnote{Obviously, the property of being 
R--nondegenerate can be  equivalently discussed  in terms of the $\Lambda$'s or in terms of the $a$'s.}

\begin{equation}\label{fmap}
\alpha : a \in \A \longmapsto  \{ \nu_1 , \dots , \nu_n , \s_1 , \dots , \s_n, \vs_1 , \dots , \vs_n    \} \in \R^{3n}
\end{equation}
where 
\begin{equation*} %\label{spazioassi}
\A :=  \{ (a_1,a_2,\dots,a_n) \in \R^n \, : \, 0 < a_n < a_{n-1} < \dots < a_1 \}  \, .
\end{equation*}

\nl
Clearly, the idea is to apply Theorem~ \ref{teoremadegenere} to the real--analytic Hamiltonian $F$ in 
(\ref{ham.n}) with $d=n$ and $p=2n$: $(\varphi,I)$ corresponding to $(\lambda,\Lambda)$ here  and  $u$ corresponding to $(\xi,q)$ and $v$ to $(\eta,p)$.
However, it turns out that {\sl the main hypothesis of 
Theorem~ \ref{teoremadegenere} does not hold}, namely, the planetary frequency map $\a$  is R--degenerate: in fact (up to rearranging the $(q,p)$--variables) one has

\begin{equation}\label{resonances}
\left\{\begin{array}{l}\vs_n=0\ , \\ \ \\
\displaystyle \sum_{j=1}^n (\s_j+\vs_j)=0\ .
\end{array}\right.
\end{equation}
The first relation is related to the rotation invariance of the system; the second relation seems to have been noticed (at least in this  generality) for the first time by Michael Herman and is therefore normally referred to as the ``Herman resonance''. 

\nl
The two resonances in (\ref{resonances}) are, however, the only linear relations identically satisfied; in fact in \cite[Proposition 78, p. 1575]{fejoz} it is proved the following

\begin{pro}
\label{prorelazionilineari}
For all $n \geq 2$ there exists an open and dense set with full Lebesgue measure $U \subset \A$,
where $\a_j\neq \a_i$ whenever $j\neq i$  and   the following property holds: 
for any open and simply connected set $V \subset U$,
the $\a_j$  define $3n$ holomorphic functions and if 
$$\a\cdot (c^1,c^2,c^3)=\nu\cdot c^1 + \s\cdot c^2+\vs \cdot c^3\equiv 0$$
for some $c^i\in\R^n$, than

\begin{equation}\label{resonances*}
\left\{\begin{array}{l}
{\rm either} \quad c^1=0\ ,\ c^2=0\ , \  c^3=(0,..,0,1)\ ,
\\ \ \\
{\rm or}\quad \quad \ \ c^1=0\ , \ c^2=(1,...,1)=c^3 \ .
\end{array}\right.
\end{equation}
\end{pro}

\nl
In order to remove the secular resonances (\ref{resonances}), 
we consider the following ``extended Hamiltonian'' on 
$$\tilde \M:=\M\times \T\times\R$$
adding a pair of conjugate symplectic variables\footnote{I.e.,   $\tilde\M$ is endowed with the symplectic form
$\displaystyle\sum_{j=1}^n \Big(d\l_j \wedge  d\Lambda_j+ 
d\xi_j \wedge  d\eta_j + dq_j \wedge  dp_j\Big) + d\theta_\rho\wedge d\rho$.
}
$(\theta_\rho,\rho)\in\T\times\R$:
\begin{equation}\label{Fmodified}
\tilde F:= F+\frac{\rho^2}2+ \e \rho^2 C_z\ \qquad{\rm with}\quad 
C_z:= \sum_{j=1}^n \Big( \Lambda_j -\frac12(\xi_j^2+\eta_j^2+q_j^2+p_j^2)\Big)\ .
\end{equation}
Let us make a few comments.

\begin{itemize}

\item[(vii)]
$C_z$ is the vertical component of the total angular momentum in Poincar\'e variables
(compare \cite{poincare} and also formula (44) in \cite{fejoz}); the form of $C_z$ is unchanged in the above variables $(\xi,\eta,q,p)$, which are obtained from the Poincar\'e variables by an orthogonal transformation.

\item[(viii)]
Since $C_z$ is an integral for $F$ (i.e., Poisson commutes with $F$), $F$ and $\tilde F$ Poisson commutes:
$$\{F,\tilde F\}\tilde{}=\{F,\tilde F \}=0
$$
where $\{\cdot,\cdot\}\tilde{}$ and $\{\cdot,\cdot\}$ denote, respectively, the Poisson bracket on $\tilde\M$ and on $\M$;
clearly, since $\tilde F$ does not depend explicitly upon the angle $\theta_\rho$, also $\rho$ is an integral for $\tilde F$ (and for $F$). 

\noindent
This fact will be  important later since 
from Lagrangian intersection theory it follows that two commuting Hamiltonians have, in general,  the same Lagrangian tori (see item (x) below for the precise statement).

\item[(ix)] The extended Hamiltonian $\tilde F$ may be rewritten as
$$\tilde F=  \tilde H^0 + \e (\tilde H^1 + H^2)
$$
with
\begin{eqnarray}
&& \tilde H^0:= F_{\rm Kep}(\Lambda)+ \frac{\rho^2}2\ ,\nonumber\\
&& \tilde H^1:= C_0(\Lambda) +  \rho^2 \sum_{j=1}^n \Lambda_j  +\\
&& \quad \quad \ \  +\sum_{j=1}^n \big(\s_j-\rho^2\big) \frac{\xi_j^2+\eta_j^2}{2} + 
														\sum_{j=1}^n \big(\varsigma_j-\rho^2\big) \frac{p_j^2+q_j^2}{2} 
+ O(4)\ .\nonumber
\end{eqnarray}
Thus, the ``slow'' action variables\footnote{Corresponding in 
Theorem~\ref{teoremadegenere}
to $I=(I_1,...,I_d)$, $d=n+1$; compare also footnote \ref{ftn*} below.} are $I=(\rho,\Lambda_1,...,\Lambda_n)$ and 
the (extended) planetary frequency map is given by
$$\tilde \a  : (\rho,a) \in \A\times \R \longmapsto  \tilde\a(\rho,a):=\Big( (\rho,\nu) ,\tilde \s, \tilde\vs\Big)\in \R^{n+1}\times\R^n\times\R^n
$$
with
\begin{equation*}
\tilde \s_j:= \s_j-\rho^2\ ,
\qquad
\tilde
\vs_j:= \vs_j-\rho^2\ . 
\end{equation*}
Proposition~\ref{prorelazionilineari}, implies immediately that $\tilde \a$ is R--nondegenerate: suppose, in fact, that 
$$\tilde \a\cdot \Big((c,c^1),c^2,c^3\Big)=  \rho\, c+\nu\cdot c^1 +  \tilde \s\cdot c^2+\tilde \vs \cdot c^3 \equiv 0$$
for some $c\in\R$ and $c^i\in\R^n$; such expression is a second order polynomial in $\rho$ and in order to vanish identically have to vanish its coefficients, i.e.,
\begin{eqnarray}
&& 
\nu\cdot c^1 + \s   \cdot c^2+ \vs \cdot c^3 =0
\ ,\label{tmp1}\\
&& 
c=0
\ ,\nonumber\\
&& -\sum_{j=1}^n c^2_j + c^3_j=0\ .\label{tmp2}
\end{eqnarray}
But then, by Proposition~\ref{prorelazionilineari} (and because of (\ref{tmp1})), one must have one of the alternatives listed in 
(\ref{resonances*}), which are incompatible with (\ref{tmp2}). 

\noindent
Thus {\sl $\tilde \a$ is R--nondegenerate} as claimed and Theorem~\ref{teoremadegenere} can be applied to the extended Hamiltonian\footnote{\label{ftn*}The correspondence with the notation of Theorem~\ref{teoremadegenere} being:
$d$
$=$
$n+1$,
$p$
$=$
$2n$,
$H_\e$
$=$
$\tilde F$,
$f$
$=$
$\tilde H^1+ H^2$,
$I$
$=$
$(\Lambda,\rho)$,
$(u,v)$
$=$
$\big((\xi,q),(\eta,p)\big)$,
$h(I)$
$=$
$H^0$,
$f_{00}$
$=$
$C_0+\rho^2 \sum_{j=1}^n \Lambda_j$,
$\omega$
$=$
$(\nu,\rho)$,
$\Omega$
$=$
$(\tilde\s,\tilde\vs)$.
} $\tilde F$, yielding, for $\e$ small enough, a positive measure set of real--analytic $(3n+1)$--dimensional Lagrangian tori in $\tilde \M$  invariant for $\tilde F$ and carrying quasi--periodic motion with Diophantine frequencies. 

\nl
The fact that $\tilde F$ is independent of $\theta_\rho$ and that $\omega_{*}:=\partial_{\rho} \tilde F=
(\rho + 4\e \rho C_z)$ is  constant along $\tilde F$--trajectories (compare point (viii) above) implies immediately that the tori ${\mathcal T}\subset \tilde \M$ obtained through Theorem~\ref{teoremadegenere} have the following parametrization
\begin{equation}\label{ftr}
{\mathcal T}:=\Big\{\big(Z(\psi,\theta_\r), \theta_\rho, \rho \big):\ (\psi,\theta_\rho)\in\T^{3n}\times\T\Big\}
\end{equation}
where $Z\in \M$ and with $\tilde F$--flow   given by
$$\phi_{\tilde F}^t \big(Z(\psi,\theta_\r), \theta_\rho, \rho \big)
=\big(Z(\psi+\o t,\theta_\r+\o_* t), \theta_\rho +\o_* t, \rho \big)\ ,
$$
for a suitable vector $\o\in\R^{3n}$, so that  $(\o,\o_*)$ forms  a Diophantine vector
in $\R^{3n+1}$.

\item[(x)] In \cite[Lemma 82, p. 1578]{fejoz} the following statement  is proved 

\noindent
{\sl If
$F$ and $G$ are two commuting Hamiltonians  and if ${\mathcal T}$ is a Lagrangian torus invariant for $F$ and with a dense $F$--orbit, then it is also $G$--invariant}.

\nl
Thus, since $\tilde F$ and $F$ (viewed as a functions on $\tilde \M$) commute, the tori
obtained in (ix) (on which any $\tilde F$--orbit is dense) are also invariant for the flow on $\tilde \M$ generated by $F$. Furthermore, the $F$--flow in $\tilde\M$ leaves {\sl both} $\theta_\rho$ and $\rho$ fixed so that, for any fixed $\theta_\r\in\T$, the $3n$--dimensional torus 
$${\mathcal T}_{\theta_\rho}:= \Big\{\big(Z(\psi,\theta_\r), \theta_\rho, \rho \big):\ \psi\in\T^{3n} \Big\}
$$
is invariant for $F$. But this means that such tori are invariant also for the $F$--flow in $\M$, finishing the proof of Theorem \ref{thm.an}. \hfill$\blacksquare$
\end{itemize}

\begin{oss}\rm 
The strategy followed here is similar to that followed in \cite{fejoz} with a few differences: first, in \cite{fejoz} $\rho$ is treated as a dumb parameter and no extended phase space is introduced  (but an extra argument is then needed to discuss the nondegeneracy of the frequency map with respect to parameters and to discuss the measure of the tori obtained); secondly, in \cite{fejoz} there is a restriction to a fixed vertical angular momentum submanifold, which is not needed here.
\end{oss}\rm

\end{document}